\documentclass[12pt, a4paper]{article}

\usepackage[utf8]{inputenc}
\usepackage[T1]{fontenc}
\usepackage[ngerman, english]{babel}
\usepackage[autostyle=true]{csquotes}

\usepackage{mlmodern}
\usepackage[nopatch=footnote]{microtype}
\usepackage[total={6in, 9.25in}]{geometry}

\usepackage{graphicx}

\usepackage[colorlinks=true, linkcolor=blue!70!black, citecolor=teal, urlcolor=blue!70!black]{hyperref}

\newcommand{\mbi}{\mathbf{i}}
\newcommand{\mbj}{\mathbf{j}}
\newcommand{\mbk}{\mathbf{k}}
\newcommand{\bbC}{\mathbb{C}}
\newcommand{\bbD}{\mathbb{D}}
\newcommand{\bbH}{\mathbb{H}}
\newcommand{\bbP}{\mathbb{P}}
\newcommand{\bbR}{\mathbb{R}}
\newcommand{\eps}{\varepsilon}

\renewcommand{\Re}{\operatorname{Re}}
\renewcommand{\Im}{\operatorname{Im}}
\newcommand{\Span}{\operatorname{span}}

\newcommand{\SE}[1][3]{\operatorname{SE}(#1)}

\usepackage{mathtools}
\usepackage{amssymb}
\usepackage{amsthm}

\theoremstyle{plain}
\newtheorem{theorem}{Theorem}[section]
\newtheorem{lemma}[theorem]{Lemma}
\newtheorem{proposition}[theorem]{Proposition}
\newtheorem{corollary}[theorem]{Corollary}

\theoremstyle{definition}
\newtheorem{definition}[theorem]{Definition}
\newtheorem{example}[theorem]{Example}

\theoremstyle{remark}
\newtheorem{remark}[theorem]{Remark}
 
\usepackage[backend=biber, sorting=nyt, maxbibnames=5, giveninits=true]{biblatex}
\usepackage{authblk}

\usepackage{orcidlink}

\title{Characterization and Kinematic Equivalence of Bennett 5-Spaces}
\author[1]{Andreas~Mair\,\orcidlink{0009-0006-9955-1655}}
\author[2,3]{Zijia~Li\,\orcidlink{0000-0002-8734-4396}}
\author[1]{Hans-Peter~Schröcker\,\orcidlink{0000-0003-2601-6695}}
\author[4]{Johannes~Siegele\,\orcidlink{0000-0002-8790-0081}}
\affil[1]{Unit of Geometry and Surveying, University of Innsbruck, 6020 Innsbruck, Austria}
\affil[2]{State Key Laboratory of Mathematical Sciences, Academy of Mathematics and Systems Science, Chinese Academy of Sciences, Beijing 100190, China}
\affil[3]{School of Mathematical Sciences, University of Chinese Academy of Sciences, Beijing 100049, China}
\affil[4]{Johann Radon Institute for Computational and Applied Mathematics (RICAM), Austrian Academy of Sciences, 4040 Linz, Austria}
\date{20 November 2026}

\begin{document}

\maketitle

\begin{abstract}
    We investigate the configuration variety of a Bennett 3R chain, a serial kinematic chain formed by three consecutive revolute (R) axes of Bennett's isogram mechanism. This variety is special in many respects; for example, it spans only a projective 5-space in the dual quaternion model of space kinematics, and the underlying 6-dimensional vector space is invariant under left and right multiplication by fixed dual quaternions. We classify all such subspaces: they form two families, each parameterized by a complex projective line. Calling two subspaces kinematically equivalent if they differ by multiplication by unit dual quaternions from the left and the right, that is, by a change of frames, we show that each subspace of one family is equivalent to one of the other and that, within one family, a single real number, possibly infinite, decides equivalence. For Bennett 3R chains, we compute this invariant in closed form from the kinematic design parameters. We further show that distinct oriented outer axes alone determine the 5-space, which is consequently shared by infinitely many chains with different middle axes. Additionally, we characterize Bennett 5-spaces by their polar lines with respect to the Study quadric and give two quadrics that cut out the configuration variety.
 \end{abstract}



\section{Introduction}\label{sec:introduction}

Among the possible models of the group of rigid body displacements \( \SE \), the Study quadric model~\cite{Selig_2005} is distinguished by combining a bilinear composition law with a minimal number of parameters~\cite{Study_1891}. It has been used to obtain a number of important results from space kinematics, mechanism science, and robotics, for example Husty's algorithm for computing the 40 solutions to the direct kinematics of Stewart--Gough platforms~\cite{Husty_1996}, the HUPF algorithm for computing the inverse kinematics of serial kinematic chains consisting of six revolute joints (6R chains)~\cite{Husty_2007}, the factorization of rational motions~\cite{Hegedus_2013a}, the computation of the 456 solutions to the synthesis of 3R chains for five prescribed poses~\cite{Hauenstein_2017}, or the proof of rationality of the singularity variety of general Stewart--Gough platforms~\cite{Coste_2020}.

The Study quadric model of \( \SE \) is embedded in the projective space \( \bbP^7(\bbR) \) over the real eight-dimensional vector space of dual quaternions. Rigid body displacements are represented by dual quaternions of non-vanishing real norm. They correspond to points of the Study quadric, a regular quadric of signature \( (4,4) \), minus the projective 3-space of dual quaternions of norm zero. Composition is just dual quaternion multiplication. Unsurprisingly, subspaces of \( \bbP^7(\bbR) \) play an important role in theoretical and applied kinematics. For example, straight lines in the Study quadric correspond to rotations about a fixed axis or translations in a fixed direction, composed with a fixed displacement. A systematic study of all subspaces contained in the Study quadric can be found in~\cite[Section~11.2]{Selig_2005}. Subspaces of dimension two are related to Bennett's isogram mechanisms~\cite{Brunnthaler_2005,Hamann_2011}, and a special 5-space characterizes line-symmetric motions~\cite{Hamann_2011,Selig_2011}. Kinematic 2R chains and chains consisting of one revolute and one prismatic (translational) joint generate particular 3-spaces that were characterized in~\cite{Rad_2018}, while~\cite{Lubbes_2018} demonstrates a kinematic generation of Darboux cyclides in terms of the general intersection of 3-spaces with the Study quadric.

The configuration variety of a 2R chain spans a projective 3-space; that of a 3R chain spans the whole space \( \bbP^7(\bbR) \)~-- in general. Exceptions are planar and spherical 3R chains, which generate 3-spaces contained in the Study quadric, and, curiously, triples of consecutive revolute axes stemming from Bennett's isogram mechanisms~\cite{Selig_2014,Hegedus_2015}. We refer to the latter as \emph{Bennett 3R chains}. They are important as building blocks of numerous overconstrained closed 6R chains, whose complete classification is one of the major open problems in theoretical mechanism science. The configuration variety of a Bennett 3R chain is the complete intersection of a 5-space with two quadrics~\cite{Selig_2014,Hegedus_2015}. These \emph{Bennett 5-spaces} are the main object of interest in this article. In future work we intend to use them for the synthesis of overconstrained closed 6R chains.

Our approach is linear-algebraic. A Bennett 5-space is mapped to itself by rotations about the first axis, applied from the left, and by rotations about the last axis, applied from the right. A change of the fixed and the moving frame can move both axes to the same coordinate axis, and then the underlying 6-dimensional real subspace of the dual quaternions is invariant under left and right multiplication by the quaternion unit \( \mbi \). In Section~\ref{sec:classification}, we classify all such subspaces: they form two families, each parameterized by a complex projective line (Theorem~\ref{thm:classification}). Since the frames can be chosen in many ways, one Bennett 5-space gives rise to many such subspaces, and we call two subspaces \( V_1, V_2 \) \emph{kinematically equivalent} if \( V_2 = pV_1 q \) for Study quaternions \( p, q \), that is, if they differ only by a change of the fixed and the moving frame. In Section~\ref{sec:equivalence}, we show that every subspace of one family is kinematically equivalent to one of the other and that, within one family, a single number \( K \in \bbR \cup \{ \infty \} \), the \emph{kinematic invariant}, decides kinematic equivalence (Theorem~\ref{thm:kinematic_equivalence}).

For Bennett 5-spaces, the invariant \( K \) is finite and has a clear connection to the kinematic design parameters (Denavit--Hartenberg parameters) of the Bennett 3R chain. In Section~\ref{sec:bennett}, we compute it from the twist angles \( \alpha_1, \alpha_2 \) and the Bennett ratio \( \ell \) of the chain, and in Section~\ref{sec:polar}, we express it by the distance \( d_{13} \) and the twist angle \( \alpha_{13} \) of the outer axes:
\begin{equation}
    \label{eq:K_intro}
    K = \tfrac{1}{2}\ell(\cos(\alpha_1) + \cos(\alpha_2)) = \tfrac{1}{2}d_{13}\cot(\tfrac{1}{2}\alpha_{13})
\end{equation}
(Theorem~\ref{thm:bennett_invariant}).
Two Bennett 5-spaces are kinematically equivalent if and only if their invariants agree. So, up to a change of frames, the 5-space reduces the three design parameters of the chain to a single number, and chains with different Denavit--Hartenberg parameters can have the same 5-space.

In Section~\ref{sec:quadrics}, we give explicit equations for two quadrics that cut out the configuration variety in the 5-space: the Study quadric, whose restriction depends only on \( K \), and a second quadric that also involves the Denavit--Hartenberg parameters (Theorem~\ref{thm:second_quadric}).

The second equality of \eqref{eq:K_intro} holds in every configuration in which the outer axes are distinct: \( d_{13} \) and \( \alpha_{13} \) change as the chain moves, but this combination of them does not. Therefore, the outer axes alone determine the 5-space; in other words, the 5-space does not depend on the middle axis (Corollary~\ref{cor:middle_axis}). The key to understanding this surprising phenomenon is the polar line of the 5-space with respect to the Study quadric. Up to the dual factor \( 1 - K\eps \), it consists of the rotations about the last axis, followed by the half-turn that interchanges the two outer axes (Theorem~\ref{thm:polar_bennett}).

\section{Preliminaries}\label{sec:preliminaries}

Let \( \bbD\bbH = \bbH \otimes_{\bbR} \bbD \) be the algebra of dual quaternions over the dual numbers \( \bbD = \bbR[\eps]/(\eps^2) \), with \( \eps^2 = 0 \). We write a dual quaternion as \( x = p + \eps d \) with \emph{primal part} \( p \in \bbH \) and \emph{dual part} \( d \in \bbH \). Let \( \bbH_{\mbi} \coloneqq \Span_{\bbR}(1, \mbi) \subset \bbH \). Because \( \mbi^2 = -1 \), the real 2-dimensional subalgebra \( \bbH_{\mbi} \) is a field isomorphic to \( \bbC \).

We use three kinds of conjugation:
\begin{itemize}
    \item For \( x = p + \eps d \in \bbD\bbH \), we write \( x^* \coloneqq p^* + \eps d^* \) for quaternion conjugation extended to \( \bbD\bbH \).
    \item For \( x \in \bbH_{\mbi} \oplus \eps\bbH_{\mbi} \), we write \( \bar{x} \) instead of \( x^* \). On \( \bbH_{\mbi} \cong \bbC \), this is the usual complex conjugation, and for \( x = x_0 + \eps x_1 \) with \( x_0, x_1 \in \bbH_{\mbi} \) it is the \( \bbD \)-linear extension \( \bar{x} = \bar{x}_0 + \eps\bar{x}_1 \). We distinguish the two notations only to make formulas that involve both more legible.
    \item The \emph{\( \eps \)-conjugate} of \( x = p + \eps d \in \bbD\bbH \) is \( x_{\eps} \coloneqq p - \eps d \).
\end{itemize}

The real part of a quaternion \( p \) is denoted by \( \Re p \coloneqq \frac{1}{2}(p+p^*) \). For a dual quaternion \( x = p + \eps d \), we call \( \Re x \coloneqq \frac{1}{2}(x + x^*) = \Re p + \eps\Re d \in \bbD \) its \emph{scalar part}. For \( q \in \bbH_{\mbi} \), we write \( q = \Re q + (\Im q)\mbi \) with \( \Im q \in \bbR \).

The quadrance (or squared norm) \( xx^* = pp^* + 2\eps\Re(pd^*) \) is in general a dual number. We also write it as \( \lvert x \rvert^2 \) and set
\begin{equation*}
    \bbD\bbH^{\times} \coloneqq \{ x \in \bbD\bbH \mid xx^* \in \bbR \setminus \{0\} \}.
\end{equation*}
The elements of \( \bbD\bbH^{\times} \) are called \emph{Study quaternions}. The condition \( xx^* \in \bbR \) is equivalent to the \emph{Study condition} \( \Re(pd^*) = 0 \), which defines the \emph{Study quadric} in \( \bbP(\bbD\bbH) \cong \bbP^7(\bbR) \). The Study quaternions \( x = p + \eps d \) correspond precisely to the points of the Study quadric minus the projective 3-space over the subspace \( \eps\bbH \) of dual quaternions with vanishing primal part.

Study quaternions serve to parameterize the group \( \SE \) of rigid body displacements. We identify a vector \( y = (y_1, y_2, y_3) \in \bbR^3 \) with the vectorial quaternion \( y_1 \mbi + y_2 \mbj + y_3 \mbk \), and the point \( y \) with the dual quaternion \( 1 + \eps y \). A Study quaternion \( x = p + \eps d \) acts on points by
\begin{equation*}
    1 + \eps y \mapsto \frac{x_{\eps}(1 + \eps y)x^*}{xx^*}.
\end{equation*}
This is the rigid body displacement represented by \( x \). Composition of displacements corresponds to the product of dual quaternions, but it can also be interpreted as a change of coordinate frames. Right multiplication changes coordinates in the pre-image space (which is called the \emph{moving frame} in kinematics); left multiplication changes coordinates in the image space (the \emph{fixed frame} in kinematics).

In the dual quaternion model, \( 1 - \frac{\eps}{2}t \) is the translation by the vector \( t \). A vectorial Study quaternion \( h \) describes a half-turn about the line with Plücker coordinates~\( h_\eps \). For a straight line with normalized Plücker coordinates \( h_\eps \) and \( \theta \in [0, 2\pi) \), the unit dual quaternion \( \cos(\frac{\theta}{2}) + \sin(\frac{\theta}{2})h \) represents the rotation by \( \theta \) about the line \( h_\eps \)~\cite{Selig_2005}. Every \( t - h \) with \( t \in \bbR \) is a non-zero scalar multiple of such a unit dual quaternion with \( \theta \neq 0 \) and hence a rotation about the same line. The case \( \theta = 0 \) is recovered in the limit \( t \to \infty \). With this understanding, we describe the set of rotations about the line \( h_\eps \) as
    \begin{equation}\label{eq:rotations_about_line}
        \{ t - h \mid t \in \bbP^1 \} = \Span_{\bbR}(1, h).
    \end{equation}
    For vectorial \( h \), we have \( h^* = -h \), and since \( h \) is
    normalized, \( h^2 = -hh^* = -1 \). In kinematic terms, this last equation
    simply states that applying a half-turn twice yields the identity.

    Normalized half-turns are used to model kinematic chains with revolute joints. In the case of an open 3R chain with oriented axes given by normalized half-turns \( h_1 \), \( h_2 \), \( h_3 \), the configuration variety is
    \begin{equation*}
        \{ (t_1 - h_1)(t_2 - h_2)(t_3 - h_3) \mid (t_1, t_2, t_3) \in {(\bbP^1)}^3 \}.
    \end{equation*}
    It is contained in the span of
    \begin{equation}\label{eq:generator}
        H \coloneqq \{ 1, h_1, h_2, h_3, h_1 h_2, h_1 h_3, h_2 h_3, h_1 h_2 h_3 \}.
    \end{equation}
    In general, this span is a vector space of dimension eight, but its dimension may drop in special cases. In any case, \( \Span(H) \) is invariant under left multiplication by \( h_1 \) and right multiplication by \( h_3 \). Indeed, taking into account \( h_1^2 = -1 \), we obtain
    \begin{multline*}
        h_1 \Span(H) = \Span( h_1, h_1^2, h_1 h_2, h_1 h_3, h_1^2 h_2, h_1^2 h_3, h_1 h_2 h_3, h_1^2 h_2 h_3 ) \\
        = \Span( h_1, -1, h_1 h_2, h_1 h_3, -h_2, -h_3, h_1 h_2 h_3, -h_2 h_3 ) = \Span(H),
    \end{multline*}
    and \( \Span(H)h_3 = \Span(H) \) follows in the same way from \( h_3^2 = -1 \).

    Furthermore, it is a classical result from space kinematics (see Study~\cite{Study_1903} or Blaschke~\cite{Blaschke_1960}, and also~\cite{Bottema_1979, Selig_2005}) that the group of Study quaternions acts transitively on oriented lines by conjugation. In particular, there exist Study quaternions \( g_1 \) and \( g_3 \) such that \( g_1 h_1 g_1^{-1} = \mbi \) and \( g_3 h_3 g_3^{-1} = \mbi \). Hence, the space \( g_1 \Span(H)g_3^{-1} \) is invariant under left and right multiplication by \( \mbi \). This is only a change of coordinates in the fixed and in the moving frame.

    \begin{definition}[Bennett 3R chain]\label{def:bennett_chain}
    Let \( L_1, L_2, L_3 \) be oriented lines in \( \bbR^3 \) such that any two consecutive lines are skew. For \( k = 1, 2 \), let \( n_k \) be the common normal of \( L_k \) and \( L_{k+1} \), oriented from \( L_{k+1} \) towards \( L_k \), let \( d_k > 0 \) be the distance between \( L_k \) and \( L_{k+1} \), and let the \emph{twist angle} \( \alpha_k \in (-\pi,\pi] \) be the angle of the right-handed rotation about \( n_k \) that carries the direction of \( L_k \) to the direction of \( L_{k+1} \). The open 3R chain with axes \( L_1, L_2, L_3 \) is a \emph{Bennett 3R chain} if the common normals \( n_1 \) and \( n_2 \) meet \( L_2 \) in the same point and the \emph{Bennett condition}
\begin{equation}\label{eq:bennett_condition}
    \ell \coloneqq \frac{d_1}{\sin(\alpha_1)} = \frac{d_2}{\sin(\alpha_2)}
\end{equation}
holds~\cite[Section~2]{Bennett_1914}. We call \( \ell \) the \emph{Bennett ratio}~\cite{Mavroidis_1995,Li_2015}.
\end{definition}

\begin{figure}[htb]
    \centering
    \includegraphics{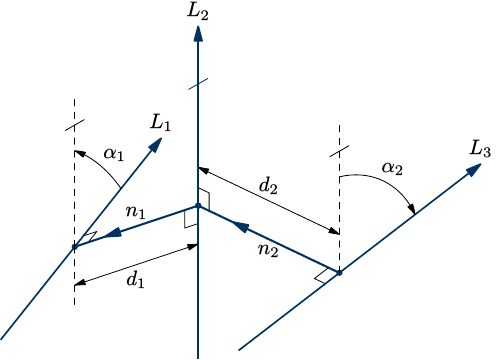}
    \caption{A Bennett 3R chain with axes \( L_1, L_2, L_3 \), common normals \( n_1, n_2 \), distances \( d_1, d_2 \) and twist angles \( \alpha_1, \alpha_2 \). Arrowheads indicate orientations. The dashed lines are parallel to \( L_2 \) and serve as references for the twist angles. Figure code written with the assistance of Claude (Anthropic).}\label{fig:chain}
\end{figure}

Figure~\ref{fig:chain} shows such a chain. In robotics, the distances \( d_k \) and the twist angles \( \alpha_k \) are known, up to notation and sign conventions, as Denavit--Hartenberg parameters~\cite{Denavit_1955}. Since consecutive axes are not parallel, we have \( \sin(\alpha_k) \neq 0 \), so \( \ell \) is well-defined. Moreover, \( \ell \neq 0 \), and \( \sin(\alpha_1), \sin(\alpha_2) \) have the sign of \( \ell \). Throughout, \( h_k \) denotes the normalized half-turn about \( L_k \), that is, the vectorial Study quaternion with \( h_k h_k^* = 1 \) such that \( h_{k,\eps} \) are the normalized Plücker coordinates of the oriented line \( L_k \). Following common usage, we also refer to \( h_1, h_2, h_3 \) as the axes of the chain.

For a Bennett 3R chain, it turns out that the eight elements of \( H \) are linearly dependent and span a 6-dimensional space~\cite{Hegedus_2013b,Selig_2014}. This is why we call the projective space over this span a Bennett 5-space and why, in the next section, we study 6-dimensional subspaces of \( \bbD\bbH \) that are invariant under left and right multiplication by~\( \mbi \). Throughout, we identify projective subspaces of \( \bbP(\bbD\bbH) \) with their underlying linear subspaces of \( \bbD\bbH \). In particular, we multiply them by Study quaternions.

\section{Classification of \texorpdfstring{\( \mbi \)}{i}-invariant subspaces}\label{sec:classification}

By Section~\ref{sec:preliminaries}, after a suitable change of the fixed and the moving frame, the underlying 6-dimensional real subspace of the Bennett 5-space of a Bennett 3R chain is invariant under left and right multiplication by \( \mbi \). In this section, we take a purely algebraic point of view and classify all such subspaces. The key is the splitting of \( \bbD\bbH \) into elements that commute and elements that anticommute with~\( \mbi \).

\begin{lemma}\label{lem:eigenspaces}
    The algebra \( \bbD\bbH \) is the direct sum \( \bbD\bbH = \mathcal{A} \oplus \mathcal{B} \) of the \( \bbH_{\mbi} \)-vector spaces
    \begin{equation*}
        \begin{aligned}
            \mathcal{A} & \coloneqq \{ q \in \bbD\bbH \mid \mbi q = q\mbi \} = \bbH_{\mbi} \oplus \eps\bbH_{\mbi},            \\
            \mathcal{B} & \coloneqq \{ q \in \bbD\bbH \mid \mbi q = -q\mbi \} = \bbH_{\mbi} \mbj \oplus \eps\bbH_{\mbi} \mbj.
        \end{aligned}
    \end{equation*}
\end{lemma}

\begin{proof}
    Any element of \( \bbH_\mbi \oplus \eps\bbH_\mbi \) commutes and any element of \( \bbH_\mbi \mbj \oplus \eps\bbH_\mbi \mbj \) anticommutes with \( \mbi \). This shows that these direct sums are contained in \( \mathcal{A} \) and \( \mathcal{B} \), respectively. Since \( \mathcal{A} \cap \mathcal{B} = \{0\} \) and the respective dimensions of the direct sums are four, equality follows.
\end{proof}

Invariant subspaces respect the splitting \( \bbD\bbH = \mathcal{A} \oplus \mathcal{B} \).

\begin{lemma}\label{lem:invariant_subspaces}
    A real subspace \( V \subset \bbD\bbH \) that is invariant under left and right multiplication by \( \mbi \) is the direct sum of \( V_{\mathcal{A}} \coloneqq V \cap \mathcal{A} \) and \( V_{\mathcal{B}} \coloneqq V \cap \mathcal{B} \).
\end{lemma}

\begin{proof}
    For \( q \in V \), also \( \mbi q \mbi \in V \), so the projections \( \pi_{\mathcal{A}}\colon q \mapsto \frac{1}{2}(q - \mbi q \mbi) \) onto \( \mathcal{A} \) and \( \pi_{\mathcal{B}}\colon q \mapsto \frac{1}{2}(q + \mbi q \mbi) \) onto \( \mathcal{B} \) map \( V \) into \( V_{\mathcal{A}} \) and \( V_{\mathcal{B}} \), respectively. As \( \pi_{\mathcal{A}} + \pi_{\mathcal{B}} = \operatorname{id} \) and \( \mathcal{A} \cap \mathcal{B} = \{0\} \), this gives \( V = V_{\mathcal{A}} \oplus V_{\mathcal{B}} \).
\end{proof}

Counting dimensions over \( \bbH_{\mbi} \) now yields the central classification result of this section.

\begin{theorem}\label{thm:classification}
    A 6-dimensional real subspace \( V \subset \bbD\bbH \) is left and right \( \mbi \)-invariant if and only if it is in one of the following two families of spaces:
    \begin{description}
        \item[Type I:] \( B_{\mathrm{I}}(u,v) \coloneqq \Span_{\bbR}( 1, \mbi, \eps, \eps\mbi, u\mbj + v\eps\mbj, u\mbk + v\eps\mbk ) \),
        \item[Type II:] \( B_{\mathrm{II}}(u,v) \coloneqq \Span_{\bbR}( u + v\eps, u\mbi + v\eps\mbi, \mbj, \mbk, \eps\mbj, \eps\mbk ) \),
    \end{description}
    where \( [u : v] \in \bbP(\bbH_{\mbi}) \) is a homogeneous parameter for each of the two families.
\end{theorem}

\begin{proof}
    Let \( V \) be a left and right \( \mbi \)-invariant subspace. By Lemma~\ref{lem:invariant_subspaces}, \( V = V_{\mathcal{A}} \oplus V_{\mathcal{B}} \) with \( \bbH_{\mbi} \)-subspaces \( V_{\mathcal{A}} \subseteq \mathcal{A} \) and \( V_{\mathcal{B}} \subseteq \mathcal{B} \). Their dimensions \( m_{\mathcal{A}}, m_{\mathcal{B}} \) over \( \bbH_{\mbi} \) satisfy \( m_{\mathcal{A}} + m_{\mathcal{B}} = 3 \), since \( \dim_{\bbR}(V) = 6 \), and \( m_{\mathcal{A}}, m_{\mathcal{B}} \le 2 \), since \( \mathcal{A} \) and \( \mathcal{B} \) are 2-dimensional over \( \bbH_{\mbi} \).

    \textit{Case 1 (\( m_{\mathcal{A}} = 2 \), \( m_{\mathcal{B}} = 1 \), Type I).} Then \( V_{\mathcal{A}} = \mathcal{A} = \Span_{\bbR}(1, \mbi, \eps, \eps\mbi) \), and \( V_{\mathcal{B}} \) is spanned over \( \bbH_{\mbi} \) by some \( b = u\mbj + v\eps\mbj \) with \( u, v \in \bbH_{\mbi} \), not both zero. A real basis of \( V_{\mathcal{B}} \) is \( b \) and \( \mbi b = u\mbk + v\eps\mbk \). Hence, \( V = B_{\mathrm{I}}(u,v) \).

    \textit{Case 2 (\( m_{\mathcal{A}} = 1 \), \( m_{\mathcal{B}} = 2 \), Type II).} Analogously, \( V_{\mathcal{B}} = \mathcal{B} = \Span_{\bbR}(\mbj, \mbk, \eps\mbj, \eps\mbk) \), and \( V_{\mathcal{A}} \) is spanned over \( \bbH_{\mbi} \) by some \( a = u + v\eps \) with \( u, v \in \bbH_{\mbi} \), not both zero. A real basis of \( V_{\mathcal{A}} \) is \( a \) and \( \mbi a = u\mbi + v\eps\mbi \). Hence, \( V = B_{\mathrm{II}}(u,v) \).

    Conversely, \( \mathcal{A} \), \( \mathcal{B} \), and the \( \bbH_{\mbi} \)-spans of \( u + v\eps \) and \( (u + v\eps)\mbj \) are invariant under left and right multiplication by \( \mbi \). Hence, so are \( B_{\mathrm{I}}(u,v) \) and \( B_{\mathrm{II}}(u,v) \). Since \( u + v\eps \neq 0 \), both are 6-dimensional.
\end{proof}

\section{Kinematic equivalence of invariant subspaces}\label{sec:equivalence}

The Study quaternions \( g_1 \) and \( g_3 \) in Section~\ref{sec:preliminaries} are not unique, so one Bennett 5-space gives rise to many \( \mbi \)-invariant subspaces. From a kinematic point of view, we should not distinguish subspaces that differ only by a change of the fixed and the moving frame. We show that every Type II subspace is equivalent in this sense to a Type I subspace and that a single number, the kinematic invariant, which may be infinite, decides equivalence within Type I.

\begin{definition}[Kinematic equivalence of subspaces]\label{def:kinematic_equivalence}
    Two real subspaces \( V_1, V_2 \subset \bbD\bbH \) are \emph{kinematically equivalent} if \( V_2 = pV_1 q \) for some \( p, q \in \bbD\bbH^\times \).
\end{definition}

\begin{lemma}\label{lem:type_exchange}
    Every Type II subspace is kinematically equivalent to a Type I subspace. More precisely, \( B_{\mathrm{I}}(u,v)\mbj = B_{\mathrm{II}}(u,v) \).
\end{lemma}

\begin{proof}
    Right multiplication by \( \mbj \) maps the generators \( 1, \mbi, \eps, \eps\mbi \) of \( B_{\mathrm{I}}(u,v) \) to \( \mbj, \mbk, \eps\mbj, \eps\mbk \) and, since \( \mbj^2 = -1 \) and \( \mbk\mbj = -\mbi \), the generators \( u\mbj + v\eps\mbj \) and \( u\mbk + v\eps\mbk \) to \( -(u + v\eps) \) and \( -(u\mbi + v\eps\mbi) \). Up to sign, these are the generators of \( B_{\mathrm{II}}(u,v) \).
\end{proof}

The remainder of this section therefore treats Type I subspaces only; the corresponding statements for Type II follow by right multiplication by \( \mbj \).

The parameter \( [u : v] \) of a Type I subspace is not preserved under kinematic equivalence. For purely imaginary \( r \in \bbH_{\mbi} \), left multiplication by the unit \( 1 + r\eps \in \mathcal{A} \) maps \( B_{\mathrm{I}}(u,v) \) to \( B_{\mathrm{I}}(u, v + ur) \) and thus replaces \( u^{-1} v \) by \( u^{-1} v + r \) if \( u \neq 0 \). Only the real part of \( u^{-1} v \) survives. This motivates the next definition.

\begin{definition}[Kinematic invariant]\label{def:kinematic_invariant}
    Let \( V = B_{\mathrm{I}}(u,v) \) be a Type I subspace with parameter \( [u : v] \in \bbP(\bbH_{\mbi}) \). The \emph{kinematic invariant} of \( V \) is
    \begin{equation*}
        K([u : v]) \coloneqq
        \begin{cases}
            \Re(u^{-1} v) & \text{if } u \neq 0, \\
            \infty        & \text{if } u = 0.
        \end{cases}
    \end{equation*}
    It is well-defined, since \( {(cu)}^{-1}(cv) = u^{-1} v \) for \( c \in \bbH_{\mbi} \setminus \{0\} \). We also write \( K(V) \coloneqq K([u : v]) \).
\end{definition}

To show that \( K \) is invariant under kinematic equivalence, we first determine which units can map one Type I subspace onto another.

\begin{lemma}\label{lem:preservation_A}
    Let \( V_1 = B_{\mathrm{I}}(u_1, v_1) \) and \( V_2 = B_{\mathrm{I}}(u_2, v_2) \), and let \( p, q \in \bbD\bbH^\times \) with \( V_2 = pV_1 q \). Then \( u_1 = 0 \) if and only if \( u_2 = 0 \). If \( u_1 \neq 0 \), then \( p \) and \( q \) either both belong to \( \mathcal{A} \) or both belong to \( \mathcal{B} \).
\end{lemma}

\begin{proof}
    If \( u_1 = 0 \), then \( V_1 \cap \eps\bbH  = \eps \bbH \). Since \( p(\eps x)q = \eps p_0 xq_0 \) for \( x \in \bbH \), where \( p_0 \) and \( q_0 \) are the primal parts of \( p \) and \( q \), respectively, the map \( x \mapsto pxq \) maps \( V_1 \cap \eps\bbH \) onto \( V_2 \cap \eps\bbH \). Comparing dimensions proves the first claim.

    Now let \( u_1 \neq 0 \), so that \( V_1 \cap \eps\bbH = \eps\bbH_{\mbi} \). Since \( V_2 \) is invariant under left multiplication by \( \mbi \), the space \( V_1 = p^{-1} V_2 q^{-1} \) is invariant under left multiplication by \( P \coloneqq p^{-1} \mbi p \). Write \( P = P_0 + \eps P_1 \) with \( P_0, P_1 \in \bbH \). Since \( 1, \eps \in V_1 \), both \( P \) and \( P \eps = \eps P_0 \) lie in \( V_1 \). Hence, \( \eps P_0 \in V_1 \cap \eps\bbH = \eps\bbH_{\mbi} \), and \( P^2 = -1 \) gives \( P_0 = \pm\mbi \). Then also \( \eps P_1 = P - P_0 \in V_1 \cap \eps\bbH \), so \( P_1 \in \bbH_{\mbi} \), and \( P^2 = -1 + 2\eps P_0 P_1 \) gives \( P_1 = 0 \). Thus, \( p^{-1} \mbi p = \pm\mbi \), that is, \( p \in \mathcal{A} \) or \( p \in \mathcal{B} \). Applying the same argument to right multiplication by \( q\mbi q^{-1} \) gives \( q \in \mathcal{A} \) or \( q \in \mathcal{B} \). If one of \( p \) and \( q \) belonged to \( \mathcal{A} \) and the other to \( \mathcal{B} \), then the 4-dimensional space \( p\mathcal{A}q \) would be contained in the 2-dimensional space \( V_2 \cap \mathcal{B} = \bbH_{\mbi}(u_2 + v_2 \eps)\mbj \).
\end{proof}

It remains to see how such \( p \) and \( q \) act on the generator \( b = (u + v\eps)\mbj \) of \( V_{\mathcal{B}} \).
We have \( bb^* = (u + v\eps)(\bar{u} + \bar{v}\eps) = {|u|}^2 + 2\Re(u\bar{v})\eps \). Since \( \Re(u\bar{v}) = {|u|}^2 \Re(u^{-1} v) \) for \( u \neq 0 \), this gives \( bb^* = 0 \) if \( u = 0 \) and
\begin{equation}\label{eq:quadrance_B}
    bb^* = {|u|}^2 \big(1 + 2K([u : v])\eps\big) \quad \text{if } u \neq 0.
\end{equation}

\begin{lemma}\label{lem:transformation_coordinates}
    Let \( b_1 = u_1 \mbj + v_1 \eps\mbj \) with \( u_1, v_1 \in \bbH_{\mbi} \), not both zero, and let \( p, q \in \bbD\bbH^\times \) both belong to \( \mathcal{A} \) or both belong to \( \mathcal{B} \). Then \( pb_1 q = u_2 \mbj + v_2 \eps\mbj \) with \( K([u_2 : v_2]) = K([u_1 : v_1]) \).
\end{lemma}

\begin{proof}
    Since \( p \) and \( q \) both belong to \( \mathcal{A} \) or both to \( \mathcal{B} \), we have \( pb_1 q \in \mathcal{B} \), that is, \( pb_1 q = u_2 \mbj + v_2 \eps\mbj \) with \( u_2, v_2 \in \bbH_{\mbi} \). As \( b_1 b_1^* \) is a dual number and \( pp^*, qq^* \in \bbR \), we have \( (pb_1 q){(pb_1 q)}^* = pb_1 qq^* b_1^* p^* = (pp^* qq^*) b_1 b_1^* \), which determines the same kinematic invariant by~\eqref{eq:quadrance_B} and homogeneity.
\end{proof}

\begin{theorem}\label{thm:kinematic_equivalence}
    Let \( V_1 \coloneqq B_{\mathrm{I}}(u_1, v_1) \) and \( V_2 \coloneqq
    B_{\mathrm{I}}(u_2, v_2) \) be two Type I \( \mbi \)-invariant subspaces.
    The subspaces \( V_1 \) and \( V_2 \) are kinematically equivalent if and
    only if \( {K([u_1 : v_1])} = {K([u_2 : v_2])} \). In this case, if \( u_1,
    u_2 \neq 0 \), then \( V_2 = pV_1 \) for the unit \( p \coloneqq 1 + (z_2 -
    z_1)\eps \in \mathcal{A} \), where \( z_k \coloneqq u_k^{-1} v_k \) for \( k
    = 1, 2\).
\end{theorem}

\begin{proof}
    Let \( b_k \coloneqq u_k \mbj + v_k \eps\mbj \), so that \( V_k = \mathcal{A} \oplus \bbH_{\mbi} b_k \) for \( k = 1, 2 \).

    First, let \( V_2 = pV_1 q \) with \( p, q \in \bbD\bbH^\times \). By Lemma~\ref{lem:preservation_A}, \( u_1 = 0 \) if and only if \( u_2 = 0 \), and in this case both invariants are \( \infty \). Otherwise, \( p \) and \( q \) both belong to \( \mathcal{A} \) or both belong to \( \mathcal{B} \), again by Lemma~\ref{lem:preservation_A}. Then \( pb_1 q \) lies in \( \mathcal{B} \) by Lemma~\ref{lem:transformation_coordinates} and in \( V_2 \), hence in \( V_2 \cap \mathcal{B} = \bbH_{\mbi} b_2 \). Thus, \( pb_1 q = cb_2 \) with \( c \in \bbH_{\mbi} \setminus \{0\} \), and Lemma~\ref{lem:transformation_coordinates} gives \( K([u_1 : v_1]) = K([cu_2 : cv_2]) = K([u_2 : v_2]) \).

    Conversely, let \( K([u_1 : v_1]) = K([u_2 : v_2]) \). If this value is \( \infty \), then \( u_1 = u_2 = 0 \) and \( V_1 = \mathcal{A} \oplus \eps\bbH_{\mbi} \mbj = V_2 \). Otherwise, \( u_1, u_2 \neq 0 \), and \( \Re(z_1) = \Re(z_2) \). Since \( z_2 - z_1 \) is purely imaginary, \( p = 1 + (z_2 - z_1)\eps \) is a unit in \( \mathcal{A} \), and
    \begin{equation*}
        pb_1 = (1 + (z_2 - z_1)\eps)u_1(1 + z_1 \eps)\mbj = u_1(1 + z_2 \eps)\mbj = u_1 u_2^{-1} b_2.
    \end{equation*}
    Hence, \( pV_1 = \mathcal{A} \oplus \bbH_{\mbi} pb_1 = \mathcal{A} \oplus \bbH_{\mbi} b_2 = V_2 \).
\end{proof}

By Theorem~\ref{thm:kinematic_equivalence} and Lemma~\ref{lem:type_exchange}, every 6-dimensional left and right \( \mbi \)-invariant subspace is kinematically equivalent to exactly one of the subspaces \( B_{\mathrm{I}}(1, K) \) with \( K \in \bbR \) or \( B_{\mathrm{I}}(0, 1) \).

\section{Bennett 5-spaces as invariant subspaces}\label{sec:bennett}

We return to Bennett 3R chains and apply the results of Sections~\ref{sec:classification} and~\ref{sec:equivalence} to their 5-spaces. After the change of the fixed and the moving frame described in Section~\ref{sec:preliminaries}, their 5-spaces always become projective spaces over Type I subspaces, and the kinematic invariant does not depend on the choice of \( g_1 \) and \( g_3 \). It is therefore an invariant of the Bennett 5-space itself. We compute it in closed form from the design parameters of the chain.

For a Bennett 3R chain with axes \( h_1, h_2, h_3 \), let \( \mathcal{S}
\coloneqq \bbP(\Span_{\bbR}(H)) \) denote its Bennett 5-space, with \( H \) as
in~\eqref{eq:generator}. We call a pair \( (g_1, g_3) \) of Study quaternions with \( g_1 h_1 g_1^{-1} = g_3 h_3 g_3^{-1} = \mbi \) \emph{admissible}. By Section~\ref{sec:preliminaries}, every admissible pair turns \( \mathcal{S} \) into an \( \mbi \)-invariant 5-space \( g_1 \mathcal{S}g_3^{-1} \). Up to kinematic equivalence, the choice of the pair does not matter.

\begin{lemma}\label{lem:conjugator_normalization}
    Let \( \mathcal{S} \) be the Bennett 5-space of a Bennett 3R chain, and let \( (g_1, g_3) \) and \( (g_1', g_3') \) be admissible pairs. Then \( g_1'\mathcal{S}{g_3'}^{-1} = p g_1 \mathcal{S}g_3^{-1} q \) with \( p, q \in \mathcal{A} \cap \bbD\bbH^\times \). In particular, the two 5-spaces are kinematically equivalent.
\end{lemma}

\begin{proof}
    Set \( p \coloneqq g_1'g_1^{-1} \) and \( q^{-1} \coloneqq g_3'g_3^{-1} \). Since \( h_1 = g_1^{-1} \mbi g_1 \), we have \( p\mbi p^{-1} = g_1'h_1{g_1'}^{-1} = \mbi \), so \( p \) commutes with \( \mbi \), that is, \( p \in \mathcal{A} \), and, likewise, \( q^{-1} \in \mathcal{A} \).
    Thus, \( g_1'\mathcal{S}{g_3'}^{-1} = p g_1 \mathcal{S}g_3^{-1} q \).
\end{proof}

Since \( p\mathcal{A}q = \mathcal{A} \) for such \( p, q \), the type of \( g_1 \mathcal{S}g_3^{-1} \) does not depend on the admissible pair, and by Theorem~\ref{thm:kinematic_equivalence}, neither does the kinematic invariant of a Type I subspace. It therefore suffices to consider a single admissible pair.

We next show that changing the frames or turning the middle joint does not change the subspaces \( g_1 \mathcal{S}g_3^{-1} \).

\begin{lemma}\label{lem:middle_joint}
    Let \( \mathcal{S} \) be the Bennett 5-space of a Bennett 3R chain with axes \( h_1, h_2, h_3 \), and let \( (g_1, g_3) \) be an admissible pair.
    \begin{enumerate}
        \renewcommand{\labelenumi}{(\roman{enumi})}
        \item A simultaneous change of the fixed and the moving frame by \( c \in \bbD\bbH^\times \), that is, conjugating the whole chain by \( c \), replaces \( \mathcal{S} \) by \( c\mathcal{S}c^{-1} \), and \( (g_1 c^{-1}, g_3 c^{-1}) \) is an admissible pair of the new chain.
        \item Turning the middle joint by a rotation \( r \in \Span_{\bbR}(1, h_2) \) replaces \( \mathcal{S} \) by \( \mathcal{S}r^{-1} \), and \( (g_1, g_3 r^{-1}) \) is an admissible pair of the new chain.
    \end{enumerate}
    In both cases, the new pair gives the same subspace \( g_1 \mathcal{S}g_3^{-1} \).
\end{lemma}

\begin{proof}
    In both cases, the new outer axes are \( ah_1 a^{-1} \) and \( bh_3 b^{-1} \), and the new 5-space is \( a\mathcal{S}b^{-1} \), where \( (a, b) = (c, c) \) in (i) and \( (a, b) = (1, r) \) in (ii). In (i), this holds because every axis is replaced by \( c h_k c^{-1} \). In (ii), only \( h_3 \) is replaced by \( r h_3 r^{-1} \), and the new products \( p_1 p_2 (r p_3 r^{-1}) = p_1 (p_2 r) p_3 r^{-1} \) with \( p_k \in \{1, h_k\} \) span \( \mathcal{S} r^{-1} \), because right multiplication by \( r \) maps \( \Span_{\bbR}(1, h_2) \) onto itself.

    Since \( g_1 a^{-1}(ah_1 a^{-1})ag_1^{-1} = g_1 h_1 g_1^{-1} = \mbi \), and likewise for \( g_3 \) and \( h_3 \), the pair \( (g_1 a^{-1}, g_3 b^{-1}) \) is admissible for the new chain, and
    \begin{equation*}
        (g_1 a^{-1})(a \mathcal{S} b^{-1}){(g_3 b^{-1})}^{-1} = g_1 \mathcal{S}g_3^{-1}. \qedhere
    \end{equation*}
\end{proof}

Lemma~\ref{lem:middle_joint} also describes how the Bennett 5-space depends on the configuration in which the axes are given.

\begin{remark}\label{rem:home_configuration}
    The axes \( h_1, h_2, h_3 \) describe the chain in one configuration, its \emph{home configuration}. It is the point \( 1 \) (the identity transformation) of the configuration variety, obtained for \( t_1 = t_2 = t_3 = \infty \), as explained right before~\eqref{eq:rotations_about_line}. The Bennett 5-space depends on this choice. To take the configuration \( x = (t_1 - h_1)(t_2 - h_2)(t_3 - h_3) \) as the new home configuration, we turn the middle joint by \( t_2 - h_2 \), thus changing the third joint, and then the first joint by \( t_1 - h_1 \), thus changing the second and third joints; turning the last joint moves no axis. By Lemma~\ref{lem:middle_joint}(ii), the first step replaces \( \mathcal{S} \) by \( \mathcal{S}{(t_2 - h_2)}^{-1} \). The second step fixes \( h_1 \) and conjugates the other axes by \( t_1 - h_1 \). It is therefore a simultaneous change of the fixed and the moving frame by \( t_1 - h_1 \), which conjugates the 5-space by \( t_1 - h_1 \) by Lemma~\ref{lem:middle_joint}(i). As \( (t_1 - h_1)\mathcal{S} = \mathcal{S} \) and \( \mathcal{S}(t_3 - h_3) = \mathcal{S} \) by Section~\ref{sec:preliminaries}, the new 5-space is
    \begin{equation*}
        (t_1 - h_1) \mathcal{S}{(t_2 - h_2)}^{-1}{(t_1 - h_1)}^{-1} = \mathcal{S}x^{-1}.
    \end{equation*}
    This is a change of the moving frame by \( x^{-1} \). In particular, the 5-spaces of all home configurations are kinematically equivalent, and by Lemma~\ref{lem:middle_joint}, the new admissible pairs give the same subspaces \( g_1 \mathcal{S}g_3^{-1} \).
\end{remark}

\begin{figure}[htb]
    \centering
    \includegraphics{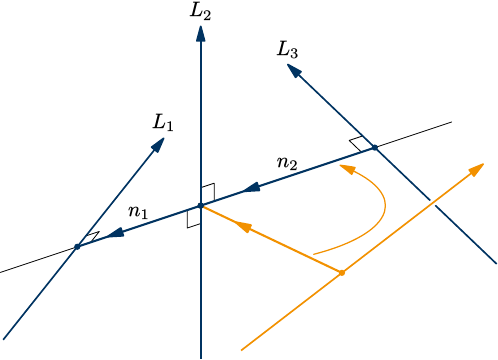}
    \caption{The chain of Figure~\ref{fig:chain} in the stretched-out configuration. The middle axis \( L_2 \) lies on the third coordinate axis, and the common normals \( n_1 \) and \( n_2 \) lie on the first. The orange lines show \( n_2 \) and \( L_3 \) in the configuration of Figure~\ref{fig:chain}, before the middle joint is turned. Figure code written with the assistance of Claude (Anthropic).}\label{fig:stretched}
\end{figure}

We now choose a configuration of the chain in which an admissible pair can be written down explicitly. By Definition~\ref{def:bennett_chain}, the common normals \( n_1 \) and \( n_2 \) meet the middle axis \( L_2 \) in the same point. A simultaneous change of the fixed and the moving frame puts \( L_2 \) on the third coordinate axis, this point at the origin, and \( n_1 \) on the first coordinate axis. Suitably turning the middle joint puts \( n_2 \) on the first coordinate axis as well. This is the ``stretched-out'' configuration of~\cite{Selig_2014}, shown in Figure~\ref{fig:stretched}. By Lemma~\ref{lem:middle_joint}, neither step changes the subspaces \( g_1 \mathcal{S}g_3^{-1} \). This configuration and the admissible pair of the next example will reappear in many later computations.

\begin{example}\label{ex:stretched}
    In the stretched-out configuration, the axes of a Bennett 3R chain are represented by the dual quaternions
    \begin{align*}
        h_1 & \coloneqq \sin(\alpha_1) \mbj + \cos(\alpha_1) \mbk + d_1 \eps ( \cos(\alpha_1) \mbj - \sin(\alpha_1) \mbk ),  \\
        h_2 & \coloneqq \mbk,                                                                                                \\
        h_3 & \coloneqq -\sin(\alpha_2) \mbj + \cos(\alpha_2) \mbk - d_2 \eps ( \cos(\alpha_2) \mbj + \sin(\alpha_2) \mbk ),
    \end{align*}
    with twist angles \( \alpha_1, \alpha_2 \) and normal distances \( d_1, d_2 \) as in Definition~\ref{def:bennett_chain}, subject to the Bennett condition~\eqref{eq:bennett_condition}. Since
    \begin{equation*}
        (\mbi + h_1) h_1 = \mbi h_1 - 1 = \mbi (\mbi + h_1)
    \end{equation*}
    and
    \begin{equation*}
        h_3 (\mbi + h_3) = h_3 \mbi - 1 = (\mbi + h_3) \mbi,
    \end{equation*}
    we have \( (\mbi + h_1) h_1 {(\mbi + h_1)}^{-1} = \mbi \) and \( {(\mbi + h_3)}^{-1} h_3 (\mbi + h_3) = \mbi \). Setting \( g_1 \coloneqq \mbi + h_1 \) and \( g_3 \coloneqq {(\mbi + h_3)}^{-1} \) therefore gives an admissible pair, provided that \( g_1 \) and \( g_3 \) are Study quaternions.

    These \( g_1 \) and \( g_3 \) are Study quaternions: both \( h_1 \) and \( h_3 \) have vanishing \( \mbi \)-component, so \( \mbi h_k + h_k \mbi = 0 \) for \( k = 1, 3 \), and since \( \mbi + h_k \) has vanishing scalar part,
    \begin{equation*}
        (\mbi + h_k){(\mbi + h_k)}^* = -{(\mbi + h_k)}^2 = 2 - (\mbi h_k + h_k \mbi) = 2.
    \end{equation*}
    Hence, \( g_1 g_1^* = 2 \) and \( g_3 g_3^* = \tfrac{1}{2} \), both real.
\end{example}

Computations with these axes simplify in complex notation. Over \( \bbH_{\mbi} \cong \bbC \), the twist angles enter only through the unit complex numbers
\begin{equation*}
    \zeta_k \coloneqq
    \cos(\alpha_k) - \sin(\alpha_k)\mbi \in \bbH_{\mbi}, \quad k = 1, 2.
\end{equation*}
Writing \( \mbk = \mbi\mbj \) and collecting the coefficients of \( \mbj \), the axes take the compact form
\begin{equation*}
    h_1 = (\mbi + d_1 \eps)\zeta_1 \mbj, \quad h_2 = \mbi\mbj, \quad h_3 = (\mbi - d_2 \eps)\bar{\zeta}_2 \mbj,
\end{equation*}
and the Bennett condition \( d_k = \ell\sin(\alpha_k) \) becomes
\begin{equation}\label{eq:bennett_compact}
    d_k = \tfrac{1}{2}\ell(\zeta_k - \bar{\zeta}_k)\mbi, \quad k = 1, 2.
\end{equation}

The next lemma records the images of two elements of \( \mathcal{S} \) under \( x \mapsto g_1 x g_3^{-1} \). It is a technical computation, but the explicit formulas are the basis for several later results.

\begin{lemma}\label{lem:stretched_images}
    For the axes and the admissible pair \( (g_1, g_3) \) of Example~\ref{ex:stretched}, the images of \( 1, h_2 \in \mathcal{S} \) under \( x \mapsto g_1 xg_3^{-1} \) are
    \begin{equation}\label{eq:stretched_images}
        \begin{aligned}
            g_1 g_3^{-1}     & = -(1 + \zeta_1 \zeta_2) + (d_1 + d_2)\zeta_1 \zeta_2 \mbi \eps + \big( \zeta_1 - \bar{\zeta}_2 - (d_1 \zeta_1 + d_2 \bar{\zeta}_2)\mbi \eps \big)\mbj,  \\
            g_1 h_2 g_3^{-1} & = -(\zeta_1 + \zeta_2)\mbi - (d_1 \zeta_1 + d_2 \zeta_2) \eps + \big( (1 - \zeta_1 \bar{\zeta}_2)\mbi - (d_1 - d_2)\zeta_1 \bar{\zeta}_2 \eps \big)\mbj.
        \end{aligned}
    \end{equation}
\end{lemma}

\begin{proof}
    We expand both products and sort the resulting terms into their \( \mathcal{A} \)- and \( \mathcal{B} \)-components. Since \( g_1 = \mbi + h_1 \) and \( g_3^{-1} = \mbi + h_3 \), expanding the products gives
    \begin{equation*}
        g_1 g_3^{-1} = -1 + \mbi h_3 + h_1 \mbi + h_1 h_3 \quad \text{and} \quad g_1 h_2 g_3^{-1} = \mbi h_2 \mbi + \mbi h_2 h_3 + h_1 h_2 \mbi + h_1 h_2 h_3.
    \end{equation*}
    We compute each term from the compact form of the axes by moving every \( \mbj \) to the right with \( \mbj c = \bar{c}\mbj \) for \( c \in \mathcal{A} \) and then using \( \mbj^2 = -1 \). For instance, \( h_1 h_3 = (\mbi + d_1 \eps)\zeta_1(-\mbi - d_2 \eps)\zeta_2 \mbj^2 \), because \( \mbj(\mbi - d_2 \eps)\bar{\zeta}_2 = (-\mbi - d_2 \eps)\zeta_2 \mbj \). This gives
    \begin{align*}
        \mbi h_3 & = -(1 + d_2 \mbi \eps)\bar{\zeta}_2 \mbj,                & \mbi h_2 \mbi & = \mbi\mbj,                                                      \\
        h_1 \mbi & = (1 - d_1 \mbi \eps)\zeta_1 \mbj,                       & \mbi h_2 h_3  & = -(\mbi + d_2 \eps)\zeta_2,                                     \\
        h_1 h_3  & = \big( -1 + (d_1 + d_2) \mbi \eps \big)\zeta_1 \zeta_2, & h_1 h_2 \mbi  & = -(\mbi + d_1 \eps)\zeta_1,                                     \\
                 &                                                          & h_1 h_2 h_3   & = -\big( \mbi + (d_1 - d_2)\eps \big)\zeta_1 \bar{\zeta}_2 \mbj.
    \end{align*}
    The terms without \( \mbj \) lie in \( \mathcal{A} \), and those with \( \mbj \) lie in \( \mathcal{B} = \mathcal{A}\mbj \). Collecting them gives~\eqref{eq:stretched_images}.
\end{proof}

By Theorem~\ref{thm:classification}, the subspace \( g_1 \mathcal{S} g_3^{-1} \) is of Type I or of Type II\@. For the admissible pair of Example~\ref{ex:stretched}, the images in Lemma~\ref{lem:stretched_images} show that it is of Type I and give its parameters explicitly.

\begin{lemma}\label{lem:stretched_type_one}
    Let \( \mathcal{S} \) be the Bennett 5-space of a Bennett 3R chain in the configuration of Example~\ref{ex:stretched}, and let \( g_1 = \mbi + h_1 \) and \( g_3 = {(\mbi + h_3)}^{-1} \). Then \( g_1 \mathcal{S}g_3^{-1} = B_{\mathrm{I}}(u,v) \) with
    \begin{equation*}
        u = \zeta_1 - \bar{\zeta}_2 \quad \text{and} \quad v = \tfrac{1}{2}\ell(\zeta_1^2 - \bar{\zeta}_2^2).
    \end{equation*}
\end{lemma}

\begin{proof}
    Let \( V \coloneqq g_1 \mathcal{S} g_3^{-1} \). By Lemma~\ref{lem:invariant_subspaces}, \( V = V_{\mathcal{A}} \oplus V_{\mathcal{B}} \), and by Theorem~\ref{thm:classification}, \( V \) is of Type I if and only if \( V_{\mathcal{A}} = \mathcal{A} \). As \( \mathcal{A} \) is 2-dimensional over \( \bbH_{\mbi} \), it suffices to find two elements of \( V \) whose \( \mathcal{A} \)-components are linearly independent over \( \bbH_{\mbi} \). We show that the images \( g_1 g_3^{-1} \) and \( g_1 h_2 g_3^{-1} \) of \( 1, h_2 \in \mathcal{S} \) have this property.

    Denote by \( a \) and \( a' \) the \( \mathcal{A} \)-components of \( g_1 g_3^{-1} \) and \( g_1 h_2 g_3^{-1} \). By~\eqref{eq:stretched_images},
    \begin{equation*}
        a = -(1 + \zeta_1 \zeta_2) + (d_1 + d_2)\zeta_1 \zeta_2 \mbi \eps \quad \text{and} \quad a' = -(\zeta_1 + \zeta_2)\mbi - (d_1 \zeta_1 + d_2 \zeta_2) \eps,
    \end{equation*}
    and the \( \mathcal{B} \)-component of \( g_1 g_3^{-1} \) is \( (u + v \eps)\mbj \) with \( u = \zeta_1 - \bar{\zeta}_2 \) and \( v = -(d_1 \zeta_1 + d_2 \bar{\zeta}_2)\mbi \).

    We write \( a = a_0 + a_1 \eps \) and \( a' = a_0' + a_1'\eps \) with \( a_0, a_1, a_0', a_1' \in \bbH_{\mbi} \). Then \( a \) and \( a' \) are linearly independent over \( \bbH_{\mbi} \) if and only if the determinant \( a_0 a_1' - a_1 a_0' \) is non-zero. Expanding and using \( 1 - \zeta_k^2 = 2\sin(\alpha_k)\mbi\zeta_k \), we obtain
    \begin{equation*}
        a_0 a_1' - a_1 a_0' = d_1 \zeta_1(1 - \zeta_2^2) + d_2 \zeta_2(1 - \zeta_1^2) = 2\mbi\zeta_1 \zeta_2 \big( d_1 \sin(\alpha_2) + d_2 \sin(\alpha_1) \big).
    \end{equation*}
    Since \( d_1, d_2 > 0 \) and \( \sin(\alpha_1) \) and \( \sin(\alpha_2) \) have the sign of \( \ell \), the bracket is non-zero. Hence, \( V_{\mathcal{A}} = \mathcal{A} \), and \( V \) is of Type I.

    It remains to determine the parameters. As \( \dim_{\bbR} V = 6 \) and \( V_{\mathcal{A}} = \mathcal{A} \), the space \( V_{\mathcal{B}} \) is 1-dimensional over \( \bbH_{\mbi} \). It contains the \( \mathcal{B} \)-component \( (u + v \eps)\mbj \) of \( g_1 g_3^{-1} \). This element is non-zero, because \( u \neq 0 \): otherwise \( \zeta_1 = \bar{\zeta}_2 \), so \( \sin(\alpha_1) = -\sin(\alpha_2) \) and hence \( d_1 = -d_2 \), contradicting \( d_1, d_2 > 0 \). Thus, \( V_{\mathcal{B}} = \bbH_{\mbi}(u + v\eps)\mbj \), that is, \( V = B_{\mathrm{I}}(u,v) \). Finally, substituting~\eqref{eq:bennett_compact} into \( v \) gives \( v = \tfrac{1}{2}\ell(\zeta_1^2 - \bar{\zeta}_2^2) \).
\end{proof}

Lemma~\ref{lem:stretched_type_one} concerns a single configuration and a single admissible pair. With the help of Lemmas~\ref{lem:conjugator_normalization} and~\ref{lem:middle_joint} we extend it to arbitrary configurations and admissible pairs.

\begin{corollary}\label{cor:always_type_one}
    Let \( \mathcal{S} \) be the Bennett 5-space of a Bennett 3R chain. For every admissible pair \( (g_1, g_3) \), the subspace \( g_1 \mathcal{S}g_3^{-1} \) is of Type I and has finite kinematic invariant.
\end{corollary}

\begin{proof}
    Let \( (g_1, g_3) \) be an admissible pair. As explained above, a simultaneous change of the frames and a turn of the middle joint bring the chain into the stretched-out configuration of Example~\ref{ex:stretched}. By Lemma~\ref{lem:middle_joint}, they turn \( (g_1, g_3) \) into an admissible pair of the new chain that gives the same subspace \( g_1 \mathcal{S}g_3^{-1} \). We may therefore assume that the chain is in the stretched-out configuration.

    Lemma~\ref{lem:stretched_type_one} determines the subspace for one particular admissible pair, namely \( V \coloneqq (\mbi + h_1)\mathcal{S}(\mbi + h_3) = B_{\mathrm{I}}(u,v) \) with \( u = \zeta_1 - \bar{\zeta}_2 \). Here \( u \neq 0 \), since \( \zeta_1 = \bar{\zeta}_2 \) would give \( d_1 = -d_2 \) by~\eqref{eq:bennett_compact}. The given pair \( (g_1, g_3) \) yields in general a different subspace, but by Lemma~\ref{lem:conjugator_normalization}, \( g_1 \mathcal{S}g_3^{-1} = pVq \) with \( p, q \in \mathcal{A} \cap \bbD\bbH^\times \). Since \( V \) contains \( \mathcal{A} \) and \( p\mathcal{A}q = \mathcal{A} \), the subspace \( pVq \) contains \( \mathcal{A} \) as well, so it is of Type I by Theorem~\ref{thm:classification}. Being kinematically equivalent to \( V \), it has the same kinematic invariant \( \Re(u^{-1} v) \) by Theorem~\ref{thm:kinematic_equivalence}, and this is finite because \( u \neq 0 \).
\end{proof}

Together with Lemma~\ref{lem:conjugator_normalization}, Corollary~\ref{cor:always_type_one} allows us to attach the kinematic invariant to the Bennett 5-space itself.

\begin{definition}[Invariant of a Bennett 5-space]\label{def:bennett_invariant}
    Let \( \mathcal{S} \) be the Bennett 5-space of a Bennett 3R chain. Its \emph{kinematic invariant} is
    \begin{equation*}
        K(\mathcal{S}) \coloneqq K\big( g_1 \mathcal{S}g_3^{-1} \big)
    \end{equation*}
    for any admissible pair \( (g_1, g_3) \). It is well-defined by Corollary~\ref{cor:always_type_one}, Lemma~\ref{lem:conjugator_normalization} and Theorem~\ref{thm:kinematic_equivalence}.
\end{definition}

The next theorem expresses \( K(\mathcal{S}) \) in terms of the design parameters of the chain.

\begin{theorem}\label{thm:bennett_invariant}
    Let \( \mathcal{S} \) be the Bennett 5-space of a Bennett 3R chain with twist angles \( \alpha_1, \alpha_2 \), normal distances \( d_1, d_2 \) and Bennett ratio \( \ell \). Then
    \begin{equation}\label{eq:bennett_invariant_value}
        K(\mathcal{S}) = \tfrac{1}{2}\ell(\cos(\alpha_1) + \cos(\alpha_2)) = \tfrac{1}{2}(d_1 \cot(\alpha_1) + d_2 \cot(\alpha_2)).
    \end{equation}
    Two Bennett 5-spaces are kinematically equivalent if and only if their invariants agree.
\end{theorem}

\begin{proof}
    By Lemma~\ref{lem:middle_joint}, the subspaces \( g_1 \mathcal{S}g_3^{-1} \), and hence \( K(\mathcal{S}) \), do not change if we place the chain in the configuration of Example~\ref{ex:stretched}; there we take \( g_1 = \mbi + h_1 \) and \( g_3 = {(\mbi + h_3)}^{-1} \). By Lemma~\ref{lem:stretched_type_one}, \( v = \frac{1}{2}\ell(\zeta_1 - \bar{\zeta}_2)(\zeta_1 + \bar{\zeta}_2) = \frac{1}{2}\ell u(\zeta_1 + \bar{\zeta}_2) \), so
    \begin{equation*}
        K(\mathcal{S}) = \Re(u^{-1} v) = \tfrac{1}{2}\ell \Re(\zeta_1 + \bar{\zeta}_2) = \frac{1}{2} \ell (\cos(\alpha_1) + \cos(\alpha_2)),
    \end{equation*}
    and the Bennett condition~\eqref{eq:bennett_condition} gives \( \ell\cos(\alpha_k) = d_k \cot(\alpha_k) \), which is the second expression in~\eqref{eq:bennett_invariant_value}. For the last statement, note that every Bennett 5-space is kinematically equivalent to its Type I subspace \( g_1 \mathcal{S}g_3^{-1} \); the claim is therefore Theorem~\ref{thm:kinematic_equivalence} applied to these subspaces.
\end{proof}

Since the twist angles and normal distances do not change when the chain moves, the kinematic invariant depends only on the geometry of the chain, not on its configuration. It also has a direct geometric meaning in terms of the outer axes.

\begin{remark}\label{rem:K_outer_axes}
    In the stretched-out configuration of Example~\ref{ex:stretched}, the rotation by \( \alpha_1 \) about the first coordinate axis followed by the translation by \( (-d_1, 0, 0) \) carries \( h_1 \) to \( h_2 \), and the rotation by \( \alpha_2 \) about the first coordinate axis followed by the translation by \( (-d_2, 0, 0) \) carries \( h_2 \) to \( h_3 \). So the first coordinate axis is also a common normal of \( h_1 \) and \( h_3 \), their distance is \( d_1 + d_2 \), and their twist angle is \( \alpha_1 + \alpha_2 \). Put
    \begin{equation*}
        \sigma \coloneqq \tfrac{1}{2}(\alpha_1 + \alpha_2), \quad \delta \coloneqq \tfrac{1}{2}(\alpha_1 - \alpha_2), \quad s \coloneqq \tfrac{1}{2}(d_1 + d_2).
    \end{equation*}
    By the sum-to-product formulas and the Bennett condition~\eqref{eq:bennett_condition}, \( d_k = \ell\sin(\alpha_k) \),
    \begin{gather*}
        \tfrac{1}{2}\ell(\cos(\alpha_1) + \cos(\alpha_2)) = \ell\cos(\sigma)\cos(\delta), \\
        s = \tfrac{1}{2}\ell(\sin(\alpha_1) + \sin(\alpha_2)) = \ell\sin(\sigma)\cos(\delta).
    \end{gather*}
    Since \( \sin(\alpha_1) \) and \( \sin(\alpha_2) \) have the sign of \( \ell \), both angles lie in \( (0,\pi) \) or both in \( (-\pi,0) \), so \( \sin(\sigma) \neq 0 \), and~\eqref{eq:bennett_invariant_value} becomes
    \begin{equation}\label{eq:bennett_invariant_value2}
        K(\mathcal{S}) = s\cot(\sigma) = \tfrac{1}{2}(d_1 + d_2)\cot\big(\tfrac{1}{2}(\alpha_1 + \alpha_2)\big).
    \end{equation}
    Thus, \( K \) is half the distance between \( h_1 \) and \( h_3 \) times the cotangent of half their twist angle.
\end{remark}

\begin{remark}\label{rem:correspondence_K_3R}
    It is a surprising insight that the kinematic invariant \( K \) determines the 5-space of a Bennett 3R chain \emph{but not the 3R chain itself}. In other words, the same 5-space, up to kinematic equivalence, belongs to different Bennett 3R chains, and different chains can have the same invariant. Theorem~\ref{thm:kinematic_equivalence} gives an explicit unit dual quaternion that carries the Type I subspace of one chain to that of the other.
\end{remark}

We illustrate Remark~\ref{rem:correspondence_K_3R} by an example.

\begin{example}\label{ex:equivalent_chains}
    Consider the two Bennett 3R chains with Denavit--Hartenberg parameters
    \begin{enumerate}
        \renewcommand{\labelenumi}{(\roman{enumi})}
        \item \( \alpha_1 = \frac{\pi}{3}, \alpha_2 = \frac{\pi}{2}, d_1 = \sqrt{3}, d_2 = 2 \),
        \item \( \alpha_1 = \frac{\pi}{4}, \alpha_2 = \frac{\pi}{2}, d_1 = 1, d_2 = \sqrt{2} \),
    \end{enumerate}
    with Bennett ratios \( \ell = 2 \) and \( \ell = \sqrt{2} \), respectively. In the configuration of Example~\ref{ex:stretched}, their axes are
    \begin{align*}
        h_1 & \coloneqq \frac{\sqrt{3}}{2} \mbj + \frac{1}{2} \mbk + \eps \left( \frac{\sqrt{3}}{2} \mbj - \frac{3}{2} \mbk \right), \\
        h_2 & \coloneqq \mbk,                                                                                                        \\
        h_3 & \coloneqq -\mbj - 2 \eps \mbk
    \end{align*}
    and
    \begin{align*}
        \tilde{h}_1 & \coloneqq \frac{\sqrt{2}}{2} \mbj + \frac{\sqrt{2}}{2} \mbk + \eps \left( \frac{\sqrt{2}}{2} \mbj - \frac{\sqrt{2}}{2} \mbk \right), \\
        \tilde{h}_2 & \coloneqq \mbk,                                                                                                                      \\
        \tilde{h}_3 & \coloneqq -\mbj - \sqrt{2} \eps \mbk.
    \end{align*}

    The chains are different, and so are their outer axes: by Remark~\ref{rem:K_outer_axes}, \( h_1, h_3 \) have distance \( 2 + \sqrt{3} \) and twist angle \( \frac{5\pi}{6} \), whereas \( \tilde{h}_1, \tilde{h}_3 \) have distance \( 1 + \sqrt{2} \) and twist angle \( \frac{3\pi}{4} \). Nevertheless,~\eqref{eq:bennett_invariant_value2} gives the same kinematic invariant
    \begin{align*}
        K & = \tfrac{1}{2}(2 + \sqrt{3})\cot\big(\tfrac{5\pi}{12}\big) = \tfrac{1}{2}(2 + \sqrt{3})(2 - \sqrt{3}) = \tfrac{1}{2}, \\
        K & = \tfrac{1}{2}(1 + \sqrt{2})\cot\big(\tfrac{3\pi}{8}\big) = \tfrac{1}{2}(1 + \sqrt{2})(\sqrt{2} - 1) = \tfrac{1}{2}.
    \end{align*}
    For the admissible pairs of Example~\ref{ex:stretched}, Lemma~\ref{lem:stretched_type_one} gives the Type I subspaces \( V_1 = B_{\mathrm{I}}(u_1, v_1) \) and \( V_2 = B_{\mathrm{I}}(u_2, v_2) \) of the two chains, where \( u^{-1} v = \frac{1}{2}\ell(\zeta_1 + \bar{\zeta}_2) \) for each chain. Hence,
    \begin{equation*}
        z_1 \coloneqq u_1^{-1} v_1 = \tfrac{1}{2} + \big(1 - \tfrac{\sqrt{3}}{2}\big)\mbi
        \quad\text{and}\quad
        z_2 \coloneqq u_2^{-1} v_2 = \tfrac{1}{2} + \tfrac{1}{2}(\sqrt{2} - 1)\mbi.
    \end{equation*}
    By Theorem~\ref{thm:kinematic_equivalence}, the unit dual quaternion
    \begin{equation*}
        p = 1 + (z_2 - z_1)\eps = 1 + \tfrac{1}{2}(\sqrt{2} + \sqrt{3} - 3)\mbi\eps
    \end{equation*}
    satisfies \( V_2 = pV_1 \). So the two 5-spaces are kinematically equivalent, although the outer axes differ; in contrast, Proposition~\ref{prop:middle_axes_exist} gives chains with equal outer axes and equal 5-spaces. More generally, by~\eqref{eq:bennett_invariant_value}, the Bennett 3R chains with a given invariant \( K \) form a two-parameter family: the three parameters \( \alpha_1, \alpha_2, \ell \) are subject to the single equation \( \frac{1}{2}\ell(\cos(\alpha_1) + \cos(\alpha_2)) = K \).
\end{example}

\section{The configuration variety of a Bennett 3R chain}\label{sec:quadrics}

By Section~\ref{sec:bennett}, the projective 5-space of a Bennett 3R chain becomes, after a change of the fixed and the moving frame, a Type I subspace \( V = B_{\mathrm{I}}(u,v) \) with finite kinematic invariant \( K \). The configuration variety of the chain then lies in \( V \), and it is the complete intersection of \( V \) with two quadrics~\cite{Selig_2014,Hegedus_2015}. These quadrics are not unique but can be taken from a pencil of quadrics. We determine two independent quadrics of this pencil in the coordinates of \( V \). The first is the Study quadric, whose restriction to \( V \) depends only on \( K \). The second necessarily depends on the remaining design parameters (Denavit--Hartenberg parameters) of the chain.

Throughout this section, \( V = B_{\mathrm{I}}(u,v) \) is a Type I subspace with \( u \neq 0 \) and kinematic invariant \( K = \Re(u^{-1} v) \). By Corollary~\ref{cor:always_type_one}, this covers the 5-spaces of Bennett 3R chains. Every \( x \in V \) can be written uniquely as
\begin{equation*}
    x = a_p + \eps a_d + c(u + v\eps)\mbj = (a_p + cu\mbj) + \eps(a_d + cv\mbj)
\end{equation*}
with \( a_p, a_d, c \in \bbH_{\mbi} \), so its primal and dual parts are \( p = a_p + cu\mbj \) and \( d = a_d + cv\mbj \).

\begin{theorem}\label{thm:study_quadric}
    The restriction of the Study quadric to \( V \) is given by
    \begin{equation*}
        Q_1(x) \coloneqq \Re(a_p \bar{a}_d) + K{|cu|}^2 = 0.
    \end{equation*}
\end{theorem}

\begin{proof}
    The Study condition is \( \Re(pd^*) = 0 \). Since \( \mbj q = \bar{q}\mbj \) for \( q \in \mathcal{A} \) and \( \mbj^2 = -1 \), we have the product rule
    \begin{equation*}
        (a + b\mbj)(a' + b'\mbj) = (aa' - b\bar{b}') + (ab' + b\bar{a}')\mbj, \quad a, a', b, b' \in \mathcal{A}.
    \end{equation*}
    As \( d^* = \bar{a}_d - cv\mbj \), it gives
    \begin{equation*}
        pd^* = (a_p + c u \mbj)(\bar{a}_d - c v \mbj) = \big( a_p \bar{a}_d + {|c|}^2 u\bar{v} \big) + (c u a_d - a_p c v) \mbj.
    \end{equation*}
    The second summand lies in \( \mathcal{B} \) and has vanishing real part. Hence, \( \Re(pd^*) = \Re(a_p \bar{a}_d) + {|c|}^2 \Re(u\bar{v}) \), and \( \Re(u\bar{v}) = {|u|}^2 \Re(u^{-1} v) = {|u|}^2 K \) gives \( \Re(pd^*) = Q_1(x) \).
\end{proof}

To find a simple second quadric, we place the chain in the stretched-out configuration of Example~\ref{ex:stretched}. This is no loss of generality: with the admissible pair adjusted as in Lemma~\ref{lem:middle_joint}, a simultaneous change of the frames and a turn of the middle joint change neither \( V \) nor the parameterization below for the chain's configuration variety, except that the turn reparameterizes \( t_2 \). We use the admissible pair \( g_1 = \mbi + h_1 \), \( g_3 = {(\mbi + h_3)}^{-1} \) of Example~\ref{ex:stretched}, for which \( V = B_{\mathrm{I}}(u,v) \) with \( u = \zeta_1 - \bar{\zeta}_2 \) by Lemma~\ref{lem:stretched_type_one}. The chain's configuration variety is parameterized by the trilinear map
\begin{equation*}
    x(t_1,t_2,t_3) = g_1(t_1-h_1)(t_2-h_2)(t_3-h_3)g_3^{-1}.
\end{equation*}
Its image \( M \), for \( (t_1,t_2,t_3) \in \bbP^1 \times \bbP^1 \times \bbP^1 \), is the configuration variety of Section~\ref{sec:preliminaries}, transformed by \( g_1 \) and \( g_3 \).

\begin{lemma}\label{lem:factorization}
    Let \( Z(t_2) \coloneqq g_1(t_2-h_2)g_3^{-1} \). Then \( x(t_1,t_2,t_3) = (t_1-\mbi) Z(t_2) (t_3-\mbi) \).
\end{lemma}

\begin{proof}
    This follows from \( g_1(t_1-h_1) = (t_1-\mbi)g_1 \) and \( (t_3-h_3)g_3^{-1} = g_3^{-1}(t_3-\mbi) \), which hold because \( g_1 h_1 g_1^{-1} = g_3 h_3 g_3^{-1} = \mbi \).
\end{proof}

Thus, the outer joints act on \( Z(t_2) \) by left and right multiplication by the rotations \( t_1 - \mbi \) and \( t_3 - \mbi \) about the first coordinate axis, and only \( Z(t_2) \) depends on the middle joint.

\begin{lemma}\label{lem:reduction}
    Write \( Z(t_2) = (P + C\mbj) + \eps(D + E\mbj) \) with \( P, C, D, E \in \bbH_{\mbi}[t_2] \). Then the coordinates of \( x(t_1,t_2,t_3) \) in \( V \) are
    \begin{equation*}
        a_p = w_1 P, \quad a_d = w_1 D, \quad cu = w_2 C,
    \end{equation*}
    where \( w_1 \coloneqq (t_1-\mbi)(t_3-\mbi) \) and \( w_2 \coloneqq (t_1-\mbi)(t_3+\mbi) \). In particular, \( {|w_1|}^2 = {|w_2|}^2 = (1+t_1^2)(1+t_3^2) \).
\end{lemma}

\begin{proof}
    Since \( \bbH_{\mbi} \) is commutative and \( \mbj q = \bar{q}\mbj \) for \( q \in \bbH_{\mbi} \), the primal part of \( x = (t_1-\mbi)Z(t_2)(t_3-\mbi) \) is
    \begin{equation*}
        p = (t_1-\mbi)(P + C\mbj)(t_3-\mbi) = (t_1-\mbi)(t_3-\mbi)P + (t_1-\mbi)\overline{(t_3-\mbi)}C\mbj = w_1 P + w_2 C\mbj.
    \end{equation*}
    Comparing with \( p = a_p + cu\mbj \) gives \( a_p = w_1 P \) and \( cu = w_2 C \). The same computation for the dual part gives \( a_d = w_1 D \).
\end{proof}

\begin{theorem}\label{thm:second_quadric}
    In the setting above, the configuration variety \( M \) is contained in the quadrics \( Q_1(x) = 0 \) and
    \begin{equation*}
        Q_2(x) \coloneqq {|a_d + \tfrac{1}{2}\ell\bar{u}a_p|}^2 - d_1^2{|cu|}^2 = 0
    \end{equation*}
    in \( V \). It is also contained in the quadric
    \begin{equation*}
        \widetilde{Q}_2(x) \coloneqq {|a_d - \tfrac{1}{2}\ell ua_p|}^2 - d_2^2{|cu|}^2 = 0.
    \end{equation*}
\end{theorem}

\begin{proof}
    Set \( \lambda \coloneqq \ell/2 \). Since \( M \) consists of rigid body displacements, \( Q_1 \) vanishes on \( M \). By Lemma~\ref{lem:reduction} and \( {|w_1|}^2 = {|w_2|}^2 = (1+t_1^2)(1+t_3^2) \),
    \begin{equation*}
        \widetilde{Q}_2 \big( x(t_1,t_2,t_3) \big) = (1+t_1^2)(1+t_3^2)\big( {|D - \lambda uP|}^2 - d_2^2{|C|}^2 \big),
    \end{equation*}
    so it suffices to show \( {|D - \lambda uP|}^2 = d_2^2{|C|}^2 \) for all \( t_2 \). Since \( Z(t_2) = t_2 g_1 g_3^{-1} - g_1 h_2 g_3^{-1} \), we read off from~\eqref{eq:stretched_images}
    \begin{align*}
        P & = -t_2(1 + \zeta_1 \zeta_2) + (\zeta_1 + \zeta_2)\mbi,              \\
        D & = t_2(d_1 + d_2)\zeta_1 \zeta_2 \mbi + (d_1 \zeta_1 + d_2 \zeta_2), \\
        C & = t_2(\zeta_1 - \bar{\zeta}_2) - (1 - \zeta_1 \bar{\zeta}_2)\mbi.
    \end{align*}
    Substituting \( d_k = \lambda(\zeta_k - \bar{\zeta}_k)\mbi \) from~\eqref{eq:bennett_compact}, a direct computation gives
    \begin{equation*}
        D - \lambda uP = \lambda(\zeta_2 - \bar{\zeta}_2)(S + T), \quad C = -\bar{\zeta}_2(S + T),
    \end{equation*}
    where \( S \coloneqq t_2(1 - \zeta_1 \zeta_2) \) and \( T \coloneqq (\zeta_2 - \zeta_1)\mbi \). As \( |\lambda(\zeta_2 - \bar{\zeta}_2)| = d_2 \) and \( |\bar{\zeta}_2| = 1 \), this proves the claim, so \( \widetilde{Q}_2 \) vanishes on \( M \).

    For \( Q_2 \), expanding the squares gives
    \begin{equation*}
        Q_2 - \widetilde{Q}_2 = 4\lambda\Re(u)\Re(a_p\bar{a}_d) - (d_1^2 - d_2^2){|cu|}^2 = 2\ell\big( \cos(\alpha_1) - \cos(\alpha_2) \big) Q_1,
    \end{equation*}
    where we used \( \Re(u) = \cos(\alpha_1) - \cos(\alpha_2) \) and, by~\eqref{eq:bennett_condition}, \( d_1^2 - d_2^2 = -\ell^2 \big( \cos(\alpha_1) - \cos(\alpha_2) \big)\big( \cos(\alpha_1) + \cos(\alpha_2) \big) \). Hence, \( Q_2 \) vanishes on \( M \) as well.
\end{proof}

The quadrics \( Q_1 \) and \( Q_2 \), and likewise \( Q_1 \) and \( \widetilde{Q}_2 \), are independent, since \( Q_2 \) and \( \widetilde{Q}_2 \) contain the term \( {|a_d|}^2 \) and \( Q_1 \) does not. Since the configuration variety is the complete intersection of \( V \) with two quadrics~\cite{Selig_2014,Hegedus_2015}, the quadrics in \( V \) containing \( M \) form a pencil. The independent quadrics \( Q_1 \) and \( Q_2 \) span this pencil, so they cut out the configuration variety in \( V \), as a complex variety. By the identity in the proof, \( Q_2 \) and \( \widetilde{Q}_2 \) agree on \( Q_1 = 0 \), so \( Q_1 = Q_2 = 0 \) and \( Q_1 = \widetilde{Q}_2 = 0 \) define the same set.

The value of Theorems~\ref{thm:study_quadric} and~\ref{thm:second_quadric} lies in the simplicity of the formulas for \( Q_1 \) and \( Q_2 \). In the coordinates of \( V \), each of them consists of two terms, and no term mixes \( a_p \) or \( a_d \) with \( c \). This reflects the action of the outer joints, which in these coordinates multiplies \( a_p \) and \( a_d \) by \( w_1 \) and \( cu \) by \( w_2 \) (Lemma~\ref{lem:reduction}). Moreover, \( Q_1 \) contains only the kinematic invariant \( K \) and is therefore the same for all Bennett 3R chains with kinematically equivalent 5-spaces. What distinguishes these chains enters only through the coefficients \( \tfrac{1}{2}\ell\bar{u} \) and \( d_1^2 \) of \( Q_2 \).

\section{The polar line of a Bennett 5-space}\label{sec:polar}

In this section we derive a simple kinematic interpretation of the polar line \( \mathcal{S}^\perp \) of a Bennett 5-space \( \mathcal{S} \) with respect to the Study quadric. It emphasizes how special Bennett 5-spaces actually are and shows that \( \mathcal{S} \) depends only on the outer axes \( h_1 \) and \( h_3 \) of the Bennett 3R chain, provided that \( h_1 \neq h_3 \). None of this refers to a special fixed or moving coordinate frame or a particular (stretched-out) configuration.

We first introduce the polar form of the Study quadric and study the polar line of a Type~I subspace. It is its \( \mathcal{B} \)-component multiplied by a dual number (Lemma~\ref{lem:polar_type_one}). For a Bennett 3R chain, this \( \mathcal{B} \)-component is spanned by \( h_1 - h_3 \) (Lemma~\ref{lem:key_identity}), and normalizing \( h_1 - h_3 \) gives a half-turn (Lemma~\ref{lem:rho_swaps}).

Polarity with respect to the Study quadric is described by the symmetric bilinear form \( \Omega\colon \bbD\bbH \times \bbD\bbH \to \bbR \) given by
\[
    (x, y) = (p_x+\eps d_x, p_y+\eps d_y) \mapsto \Re(p_x d_y^*+p_y d_x^*).
\]
For \( x = p + \eps d \), we have \( \Omega(x,x) = 2\Re(pd^*) \), so the quadrance formula of Section~\ref{sec:preliminaries} reads
\begin{equation*}
    xx^* = {|p|}^2 + \eps \Omega(x,x),
\end{equation*}
whence the Study quadric is the zero set of \( \Omega(x,x) \). As \( \Re(pq^*) \) is the Euclidean inner product of \( p, q \in \bbH \cong \bbR^4 \), the Gram matrix of \( \Omega \) in the coordinates \( (p, d) \) is \( \left[\begin{smallmatrix} 0 & I_4 \\ I_4 & 0 \end{smallmatrix}\right] \) and \( \Omega \) is non-degenerate.

For a real subspace \( U \subseteq \bbD\bbH \), we write
\[
    U^{\perp} \coloneqq \{ y \in \bbD\bbH \mid \Omega(x,y) = 0 \text{ for all } x \in U \}.
\]
Every real subspace \( U \subseteq \bbD\bbH \) satisfies \( \dim_{\bbR} U^{\perp} = 8 - \dim_{\bbR} U \) and \( {(U^{\perp})}^{\perp} = U \). In particular, for a 6-dimensional subspace \( V \), the space \( V^{\perp} \) is 2-dimensional. We call the line \( \bbP(V^{\perp}) \subset \bbP^7 \) the \emph{polar line} of the 5-space \( \bbP(V) \). It determines the 5-space, since \( V = {(V^{\perp})}^{\perp} \).

Multiplication by Study quaternions \( q_1 \) and \( q_2 \) from the left and the right, and in particular a change of the fixed and the moving frame, preserves \( \Omega \) up to a real factor. Consequently, \( {(q_1 Uq_2)}^{\perp} = q_1 U^{\perp} q_2 \) for every real subspace \( U \subseteq \bbD\bbH \). In particular, the polar line of a Bennett 5-space \( \mathcal{S} \) is obtained from that of its Type I subspace \( g_1 \mathcal{S}g_3^{-1} \) with admissible \( g_1 \) and \( g_3 \). We therefore start our considerations with Type I subspaces.

\begin{lemma}\label{lem:polar_type_one}
    Let \( V = B_{\mathrm{I}}(u,v) \) be a Type I subspace with \( u \neq 0 \), and let \( K = K(V) \) be its kinematic invariant. Then \( V^{\perp} = (1 - 2K\eps) V_{\mathcal{B}} \).
\end{lemma}

\begin{proof}
    Both spaces, \( V^\perp \) and \( V_{\mathcal{B}} \), are 2-dimensional, so it suffices to show that \( \Omega(x, y) = 0 \) for all \( x \in V \) and \( y \in (1 - 2K\eps)V_{\mathcal{B}} \). By Theorem~\ref{thm:classification} and Definition~\ref{def:kinematic_invariant}, \( V = \mathcal{A} \oplus \bbH_{\mbi} b \) with \( b \coloneqq (u + v\eps)\mbj \), and by~\eqref{eq:quadrance_B}, \( bb^* = {|u|}^2(1 + 2K\eps) \). So \( y = (1 - 2K\eps)\mu b \) with \( \mu \in \bbH_{\mbi} \). Recall that \( \Omega(x, y) \) is the dual part of \( \Re(xy^*) \). We now have:
    \begin{itemize}
        \item For \( x = \lambda b \) with \( \lambda \in \bbH_{\mbi} \), we have
              \begin{equation*}
                  xy^* = (1 - 2K\eps) \lambda bb^* \bar{\mu} = {|u|}^2(1 - 2K\eps)(1 + 2K\eps) \lambda\bar{\mu} = {|u|}^2 \lambda\bar{\mu}
              \end{equation*}
              because \( bb^* \) is central. This lies in \( \bbH_{\mbi} \), so its scalar part is real.
        \item For \( x \in \mathcal{A} \), the product \( xy^* \) lies in \( \mathcal{B} \), since \( y^* \) does, and therefore has no scalar part.
    \end{itemize}
    In both cases, \( \Omega(x, y) = 0 \).
\end{proof}

We now return to Bennett 3R chains and fix, for the rest of this section, one with axes \( h_1, h_2, h_3 \), where \( h_1 \neq h_3 \). We denote by \( \mathcal{S} \) its Bennett 5-space, by \( K = K(\mathcal{S}) \) its invariant, and by \( v_1, v_3 \) the directions of \( h_1, h_3 \), that is, their primal parts. Since \( h_1^2 = h_3^2 = -1 \), the difference \( X \coloneqq h_1 - h_3 \) satisfies
\begin{equation*}
    Xh_3 = h_1 h_3 + 1 = -h_1 X.
\end{equation*}

In Section~\ref{sec:bennett}, we computed \( K \) in the stretched-out configuration, which is reached by turning the middle joint. This changes the relative position of \( h_1 \) and \( h_3 \), which we would like to preserve here. We therefore do not turn the middle joint and only conjugate. By Lemma~\ref{lem:middle_joint}, a simultaneous change of the fixed and the moving frame by \( c \in \bbD\bbH^\times \) replaces each \( h_k \) by \( c h_k c^{-1} \) and \( \mathcal{S} \) by \( c\mathcal{S}c^{-1} \), and hence \( \mathcal{S}^\perp \) by \( c\mathcal{S}^\perp c^{-1} \). It preserves \( K \) by Definition~\ref{def:bennett_invariant}.

We may work in any fixed frame, and we choose one whose first coordinate axis is a common normal of \( h_1 \) and \( h_3 \). With this choice, their primal and dual parts lie in \( \Span_{\bbR}(\mbj, \mbk) \), that is, \( h_1, h_3 \in \mathcal{B} \). Since the argument of Example~\ref{ex:stretched} uses only this and \( h_1^2 = h_3^2 = -1 \), the elements \( g_1 \coloneqq \mbi + h_1 \) and \( g_3 \coloneqq {(\mbi + h_3)}^{-1} \) are Study quaternions and form an admissible pair. By Corollary~\ref{cor:always_type_one} and Definition~\ref{def:bennett_invariant}, \( V \coloneqq g_1 \mathcal{S}g_3^{-1} \) is a Type I subspace \( B_{\mathrm{I}}(u,v) \) with \( u \neq 0 \) and \( K(V) = K \).

By Lemma~\ref{lem:polar_type_one}, the polar line of \( V \) is determined by its \( \mathcal{B} \)-component and its kinematic invariant. This component is spanned by \( X \), and \( g_1 \) and \( g_3 \) preserve \( X \) up to a real factor:

\begin{lemma}\label{lem:key_identity}
    The \( \mathcal{B} \)-component of \( V \) is \( V_{\mathcal{B}} = \bbH_{\mbi} X \), and
    \begin{equation*}
        g_1 Xg_3^{-1} = 2X, \quad \text{equivalently} \quad g_1^{-1} Xg_3 = \tfrac{1}{2}X.
    \end{equation*}
\end{lemma}

\begin{proof}
    We find \( X \) in the \( \mathcal{B} \)-component of the image \( g_1 g_3^{-1} \) of \( 1 \in \mathcal{S} \), and then compute \( g_1 Xg_3^{-1} \) directly. Throughout, we use that \( h_1 \), \( h_3 \) and \( X \) lie in \( \mathcal{B} \) and anticommute with \( \mbi \), while products of two elements of \( \mathcal{B} \), such as \( h_1 h_3 \), lie in \( \mathcal{A} \) and commute with~\( \mbi \).

    Since \( g_1 = \mbi + h_1 \) and \( g_3^{-1} = \mbi + h_3 \), expanding gives
    \begin{equation*}
        g_1 g_3^{-1} = (-1 + h_1 h_3) + (\mbi h_3 + h_1 \mbi) = (-1 + h_1 h_3) + X\mbi
    \end{equation*}
    with \( -1 + h_1 h_3 \in \mathcal{A} \) and \( X\mbi \in \mathcal{B} \). As \( g_1 g_3^{-1} \in V \), Lemma~\ref{lem:invariant_subspaces} gives \( X\mbi \in V_{\mathcal{B}} \), and as \( V_{\mathcal{B}} \) is an \( \bbH_{\mbi} \)-subspace, also \( X = \mbi(X\mbi) \in V_{\mathcal{B}} \). Since \( h_1 \neq h_3 \), we have \( X \neq 0 \). As \( V \) is of Type I, \( V_{\mathcal{B}} \) is 1-dimensional over \( \bbH_{\mbi} \), so \( V_{\mathcal{B}} = \bbH_{\mbi} X \).

    Similarly,
    \begin{equation*}
        g_1 Xg_3^{-1} = (\mbi + h_1)X(\mbi + h_3) = \mbi X\mbi + (\mbi Xh_3 + h_1 X\mbi) + h_1 Xh_3.
    \end{equation*}
    The first term is \( \mbi X\mbi = -X\mbi^2 = X \). The element \( Xh_3 = 1 + h_1 h_3 = -h_1 X \) lies in \( \mathcal{A} \) and commutes with \( \mbi \). Hence, the bracket is \( \mbi Xh_3 - Xh_3 \mbi = 0 \), and the last term is \( h_1 Xh_3 = -Xh_3^2 = X \). This gives \( g_1 Xg_3^{-1} = 2X \), and multiplying by \( g_1^{-1} \) from the left and by \( g_3 \) from the right gives \( g_1^{-1} Xg_3 = \tfrac{1}{2}X \).
\end{proof}

Together with the quadrance formula~\eqref{eq:quadrance_B}, Lemma~\ref{lem:key_identity} shows that the outer axes alone determine~\( K \).

\begin{corollary}\label{cor:outer_axes}
    The directions of the outer axes satisfy \( v_1 \neq v_3 \), and
    \begin{equation*}
        XX^* = {|v_1 - v_3|}^2(1 + 2K\eps).
    \end{equation*}
\end{corollary}

\begin{proof}
    By Lemma~\ref{lem:key_identity} and Definition~\ref{def:kinematic_invariant}, \( X \) is a non-zero element of \( V_{\mathcal{B}} = \bbH_{\mbi}(u + v\eps)\mbj \), so \( X = (\lambda u + \lambda v\eps)\mbj \) with \( \lambda \in \bbH_{\mbi} \setminus \{0\} \). As \( \lambda u \neq 0 \) and \( K([\lambda u : \lambda v]) = K(V) = K \), formula~\eqref{eq:quadrance_B} gives \( XX^* = {|\lambda u|}^2(1 + 2K\eps) \). The primal part of \( X \) is \( v_1 - v_3 = \lambda u\mbj \), whose norm is \( |\lambda u| \neq 0 \). Hence, \( v_1 \neq v_3 \) and \( |\lambda u| = |v_1 - v_3| \).
\end{proof}

By Corollary~\ref{cor:outer_axes}, the vector \( v_1 - v_3 \) is non-zero and can serve as a direction.

\begin{definition}[Bisecting axis]\label{def:bisecting_axis}
    The \emph{bisecting axis} of the chain is the line with direction \( v_1 - v_3 \) through the midpoint of the feet of a common normal of \( h_1 \) and \( h_3 \).
\end{definition}

It meets this common normal perpendicularly, since \( v_1 \) and \( v_3 \) are perpendicular to it. As \( v_1 \) and \( v_3 \) are unit vectors, the direction \( v_1 - v_3 \) bisects the angle between \( v_1 \) and \( -v_3 \), hence the name. If \( h_1 \) and \( h_3 \) are parallel, then \( v_3 = -v_1 \) by Corollary~\ref{cor:outer_axes}, and the midpoints of the feet of all their common normals form a line with direction \( v_1 = \tfrac{1}{2}(v_1 - v_3) \). So the bisecting axis does not depend on the choice of the common normal, which is unique otherwise.

Normalizing \( X \) gives the half-turn about the bisecting axis.

\begin{lemma}\label{lem:rho_swaps}
    The dual quaternion
    \begin{equation*}
        \rho \coloneqq \frac{1 - K\eps}{|v_1 - v_3|} X
    \end{equation*}
    is the normalized half-turn about the bisecting axis with direction \( (v_1 - v_3)/|v_1 - v_3| \). It satisfies \( \rho h_3 \rho^{-1} = -h_1 \) and \( \rho h_1 \rho^{-1} = -h_3 \), that is, it interchanges \( h_1 \) and \( h_3 \), reversing their orientations.
\end{lemma}

\begin{proof}
    We first show that \( \rho \) is a normalized half-turn and then derive the two identities, which also locate its axis.

    The factor \( (1 - K\eps)/|v_1 - v_3| \) is a dual number, so it is central and equal to its conjugate. By Corollary~\ref{cor:outer_axes},
    \begin{equation*}
        \rho\rho^* = \frac{{(1 - K\eps)}^2}{{|v_1 - v_3|}^2} XX^* = {(1 - K\eps)}^2(1 + 2K\eps) = 1.
    \end{equation*}
    Like \( h_1 \) and \( h_3 \), the element \( \rho \) has vanishing scalar part. By Section~\ref{sec:preliminaries}, it is therefore a normalized half-turn, and its direction is its primal part \( (v_1 - v_3)/|v_1 - v_3| \).

    Multiplying \( Xh_3 = -h_1 X \) by the central factor gives \( \rho h_3 = -h_1 \rho \), that is, \( \rho h_3 \rho^{-1} = -h_1 \). As \( \rho^2 = -\rho\rho^* = -1 \) is central, conjugating twice by \( \rho \) is the identity map. So conjugating \( \rho h_3 \rho^{-1} = -h_1 \) by \( \rho \) gives \( h_3 = -\rho h_1 \rho^{-1} \).

    By Section~\ref{sec:preliminaries}, conjugation by \( \rho \) maps the half-turn about a line to the half-turn about the image of that line, and \( -h \) is the half-turn about the line of \( h \) with reversed orientation. So the two identities show that \( \rho \) interchanges \( h_1 \) and \( h_3 \), reversing their orientations. Hence, it also interchanges the feet of a common normal and fixes their midpoint, so its axis is the bisecting axis.
\end{proof}

We can now determine the polar line.

\begin{theorem}\label{thm:polar_bennett}
    Let \( \mathcal{S} \) be the Bennett 5-space of a Bennett 3R chain with axes \( h_1, h_2, h_3 \), where \( h_1 \neq h_3 \), let \( K = K(\mathcal{S}) \), and let \( \rho \) be the normalized half-turn about the bisecting axis of \( h_1 \) and \( h_3 \). Then
    \begin{equation}\label{eq:polar_line_formula}
        \begin{aligned}
            \mathcal{S}^{\perp} & = (1 - 2K\eps)(h_1 - h_3)\Span_{\bbR}(1, h_3) = (1 - 2K\eps)\Span_{\bbR}(1, h_1)(h_1 - h_3) \\
                                & = (1 - K\eps) \rho\Span_{\bbR}(1, h_3) = (1 - K\eps)\Span_{\bbR}(1, h_1) \rho.
        \end{aligned}
    \end{equation}
    Hence, up to the dual factor \( 1 - K\eps \), the points of the polar line are the rotations about \( h_3 \) followed by the half-turn \( \rho \), or equivalently \( \rho \) followed by the rotations about \( h_1 \).
\end{theorem}

\begin{proof}
    We compute the polar line of \( V \) and carry it back to \( \mathcal{S} \) by left-multiplying by \( g_1^{-1} \) and right-multiplying by \( g_3 \), using the frame and the notation fixed before Lemma~\ref{lem:key_identity}. The normalized half-turn about the bisecting axis is given by \( \rho \) as in Lemma~\ref{lem:rho_swaps}.

    By Lemmas~\ref{lem:polar_type_one} and~\ref{lem:key_identity}, \( V^{\perp} = (1 - 2K\eps) \bbH_{\mbi} X \), and \( \bbH_{\mbi} X = X \bbH_{\mbi} \), since \( X \in \mathcal{B} \) anticommutes with \( \mbi \). As \( \mathcal{S} = g_1^{-1} Vg_3 \), we have \( \mathcal{S}^{\perp} = g_1^{-1} V^{\perp} g_3 \), as explained before Lemma~\ref{lem:polar_type_one}. From \( g_3 h_3 g_3^{-1} = \mbi \), we get \( \bbH_{\mbi} g_3 = g_3 \Span_{\bbR}(1, h_3) \), and Lemma~\ref{lem:key_identity} gives \( g_1^{-1} Xg_3 = \frac{1}{2}X \). Hence,
    \begin{equation*}
        \mathcal{S}^{\perp} = (1 - 2K\eps) g_1^{-1} X \bbH_{\mbi} g_3 = \tfrac{1}{2}(1 - 2K\eps)X\Span_{\bbR}(1, h_3).
    \end{equation*}
    The real span absorbs the factor \( \frac{1}{2} \), which gives the first expression. As \( Xh_3 = -h_1 X \), we have \( X\Span_{\bbR}(1, h_3) = \Span_{\bbR}(1, h_1) X \), which gives the second.

    Finally, \( X = |v_1 - v_3|(1 + K\eps)\rho \) by the definition of \( \rho \), and \( (1 - 2K\eps)(1 + K\eps) = 1 - K\eps \). Substituting \( X \) into the first two expressions and absorbing the real factor \( |v_1 - v_3| \) gives the expressions with \( \rho \). The last statement follows from Section~\ref{sec:preliminaries}: the non-zero elements of \( \Span_{\bbR}(1, h) \) are the rotations about \( h \), and a product \( xx' \) is the displacement \( x' \) followed by \( x \).
\end{proof}

In particular, the middle axis does not affect the Bennett 5-space. The hypothesis \( h_1 \neq h_3 \) is needed: the proof divides by \( X = h_1 - h_3 \), and for \( h_1 = h_3 \) the 5-space does change with the middle axis.

\begin{corollary}\label{cor:middle_axis}
    Two Bennett 3R chains with the same outer axes \( h_1 \neq h_3 \) have the same Bennett 5-space.
\end{corollary}

\begin{proof}
    By Corollary~\ref{cor:outer_axes}, applied to each chain, the two chains have the same invariant \( K \), since \( v_1 \neq v_3 \) and the dual part of \( XX^* \) is \( 2K{|v_1 - v_3|}^2 \). By Theorem~\ref{thm:polar_bennett}, they therefore have the same polar line \( \mathcal{S}^{\perp} \), and hence the same Bennett 5-space \( \mathcal{S} = {(\mathcal{S}^{\perp})}^{\perp} \).
\end{proof}

Corollary~\ref{cor:middle_axis} is not vacuous: given the outer axes, there are infinitely many possible middle axes.

\begin{proposition}\label{prop:middle_axes_exist}
    Let \( h_1, h_3 \) be oriented lines with a common normal of length \( d > 0 \) and twist angle \( \alpha \in (-\pi, \pi] \setminus \{0\} \), as in Remark~\ref{rem:K_any_configuration}. Then there are infinitely many Bennett 3R chains with outer axes \( h_1, h_3 \) and pairwise distinct middle axes. By Corollary~\ref{cor:middle_axis}, they have the same Bennett 5-space, and its kinematic invariant is \( K = \frac{1}{2}d\cot(\frac{1}{2}\alpha) \).
\end{proposition}

\begin{proof}
    Put \( \sigma \coloneqq \frac{1}{2}\alpha \) and, for \( \delta \in (-\lvert\sigma\rvert, \lvert\sigma\rvert) \),
    \begin{equation*}
        \alpha_1 \coloneqq \sigma + \delta, \quad \alpha_2 \coloneqq \sigma - \delta, \quad \ell \coloneqq \frac{d}{2\sin(\sigma)\cos(\delta)}, \quad d_k \coloneqq \ell\sin(\alpha_k).
    \end{equation*}
    The angles \( \alpha_k \in (-\pi, \pi) \) have the sign of \( \sigma \), so \( \sin(\alpha_k) \) and \( \ell \) have the same sign, \( d_k > 0 \), and \( d_1/\sin(\alpha_1) = d_2/\sin(\alpha_2) = \ell \). These are the parameters of a Bennett 3R chain. In its stretched-out configuration, the outer axes have distance \( d_1 + d_2 = 2\ell\sin(\sigma)\cos(\delta) = d \) and twist angle \( \alpha_1 + \alpha_2 = \alpha \) by Remark~\ref{rem:K_outer_axes}, independently of \( \delta \). A pair of oriented lines is determined up to a displacement by these two numbers, so after a displacement the outer axes are \( h_1, h_3 \). Different \( \delta \) give different twist angles \( \alpha_1 \) between \( h_1 \) and the middle axis, which differ by less than \( \pi \), hence different middle axes.

    Finally, \( \alpha \neq 0 \) gives \( h_1 \neq h_3 \), so Corollary~\ref{cor:middle_axis} applies, and \( K \) is given by~\eqref{eq:K13}.
\end{proof}

Note that this constructive proof does not yield all possible middle axes.

\begin{proposition}\label{prop:polar_characterization}
    A 6-dimensional real subspace \( V \subset \bbD\bbH \) is kinematically equivalent to a Bennett 5-space with kinematic invariant \( K \in \bbR \) if and only if
    \begin{equation*}
        V^{\perp} = (1 - K\eps)\Span_{\bbR}(1, h)\,g
    \end{equation*}
    for a normalized half-turn \( h \) and a Study quaternion \( g \).
\end{proposition}

\begin{proof}
    By Lemma~\ref{lem:polar_type_one}, \( {B_{\mathrm{I}}(1, K)}^{\perp} = (1 - 2K\eps)\bbH_{\mbi}(1 + K\eps)\mbj = (1 - K\eps)\bbH_{\mbi}\mbj \). This has the stated form, and multiplication by Study quaternions from the left and the right preserves both this form and the factor \( 1 - K\eps \). By Corollary~\ref{cor:always_type_one} and Theorem~\ref{thm:kinematic_equivalence}, a Bennett 5-space with invariant \( K \) is kinematically equivalent to \( B_{\mathrm{I}}(1, K) \), which proves the first direction.

    Conversely, choose \( c \in \bbD\bbH^\times \) with \( chc^{-1} = \mbi \), and set \( q \coloneqq {(cg)}^{-1}\mbj \). Then \( cV^{\perp}q = (1 - K\eps)\bbH_{\mbi}\mbj \), so \( cVq = B_{\mathrm{I}}(1, K) \). By Theorem~\ref{thm:bennett_invariant}, every real \( K \) is the invariant of some Bennett 3R chain, whose 5-space is therefore kinematically equivalent to \( V \).
\end{proof}

\begin{remark}\label{rem:K_any_configuration}
    The description of \( K \) at the end of Remark~\ref{rem:K_outer_axes} holds in every configuration of the chain, provided \( h_1 \neq h_3 \). Let \( n \) be a common normal of \( h_1 \) and \( h_3 \) with feet \( F_1 \) and \( F_3 \), oriented from \( F_3 \) towards \( F_1 \) if \( F_1 \neq F_3 \). As in Definition~\ref{def:bennett_chain}, let \( d_{13} \ge 0 \) be the distance between \( F_1 \) and \( F_3 \), and let \( \alpha_{13} \in (-\pi, \pi] \) be the angle of the right-handed rotation about \( n \) that carries \( v_1 \) to \( v_3 \), which is non-zero by Corollary~\ref{cor:outer_axes}. With the Plücker coordinates of Section~\ref{sec:preliminaries}, \( h_k = v_k + \eps(v_k \times F_k) \). Using that the scalar part of the product of two vectorial quaternions \( p, q \) is \( -p \cdot q \), we obtain
    \begin{align*}
        XX^* & = 2 + h_1 h_3 + h_3 h_1 = 2 - 2v_1 \cdot v_3 + 2\eps(F_1 - F_3) \cdot (v_1 \times v_3) \\
             & = 4\sin^2(\tfrac{1}{2}\alpha_{13}) + 2\eps d_{13} \sin(\alpha_{13}),
    \end{align*}
    since \( v_1 \cdot v_3 = \cos(\alpha_{13}) \), \( F_1 - F_3 = d_{13} n \) and \( v_1 \times v_3 = \sin(\alpha_{13})n \). Comparing with Corollary~\ref{cor:outer_axes} gives
    \begin{equation}\label{eq:K13}
        K = \tfrac{1}{2}d_{13} \cot(\tfrac{1}{2}\alpha_{13}).
    \end{equation}
    When the middle joint turns, \( d_{13} \) and \( \alpha_{13} \) change, but by Lemma~\ref{lem:middle_joint}, \( K \) does not. Moreover, \( K = 0 \) if and only if \( d_{13} = 0 \) or \( \alpha_{13} = \pi \), that is, if and only if \( h_1 \) and \( h_3 \) are coplanar. By Theorem~\ref{thm:bennett_invariant}, this happens exactly for \( |\alpha_1 + \alpha_2| = \pi \), and then the outer axes are coplanar in every configuration.
\end{remark}

\begin{remark}\label{rem:polar_K}
    By Theorem~\ref{thm:polar_bennett}, every non-zero \( x \in \mathcal{S}^{\perp} \) is \( x = (1 - K\eps)G \), where the Study quaternion \( G \) is a rotation about \( h_3 \) followed by \( \rho \). Its quadrance \( xx^* = (1 - 2K\eps) GG^* \) has the positive primal part \( GG^* \) and the dual part \( -2K GG^* \). Hence, the polar line of \( \mathcal{S} \) lies in the Study quadric if \( K = 0 \), that is, if \( h_1 \) and \( h_3 \) are coplanar (Remark~\ref{rem:K_any_configuration}), and it is disjoint from the Study quadric otherwise. In the first case, its points are the rotations about \( h_3 \) followed by the half-turn \( \rho \). In general, comparing with \( xx^* = {|p_x|}^2 + \eps \Omega(x,x) \), where \( p_x \) is the primal part of \( x \), shows that
    \begin{equation*}
        K = -\frac{\Omega(x,x)}{2{|p_x|}^2}
    \end{equation*}
    can be read off any point \( x \) of the polar line.
\end{remark}

\begin{remark}\label{rem:darboux}
    In the interpretation of space kinematics within the projective geometry of \( \bbD\bbH \) over the dual numbers, a straight line in \( \bbP(\bbD\bbH) \) describes a \emph{vertical Darboux motion}, that is, a constant-speed rotation combined with a harmonic oscillation~\cite{Siegele_2021}. The \emph{back-projection} of a point \( x \) with non-vanishing primal part is the unique point of the Study quadric of the form \( (1 + c\eps)x \) with \( c \in \bbR \). The back-projection of a straight line is, in general, a twisted cubic, but it reduces to a straight line if the amplitude of the oscillation vanishes, that is, if the vertical Darboux motion is just a rotation.

    This is the case for the polar line of a Bennett 5-space. By Remark~\ref{rem:polar_K}, it is disjoint from the Study quadric unless \( K = 0 \), but by Theorem~\ref{thm:polar_bennett}, back-projection removes the dual factor \( 1 - K\eps \) and yields the straight line \( \rho\Span_{\bbR}(1, h_3) \). Conversely, if the polar line \( V^\perp \) of a 5-space back-projects to a line \( \Span_{\bbR}(1, h)g \), then the dual factors of the back-projection agree for all its points, since \( V^\perp \) is linear and the primal parts of \( \Span_{\bbR}(1, h)g \) are linearly independent. Hence, by Proposition~\ref{prop:polar_characterization}, Bennett 5-spaces are, up to kinematic equivalence, the only 6-dimensional subspaces whose polar line back-projects to rotations about a fixed axis composed with a fixed displacement.
\end{remark}

\section{Conclusion}\label{sec:conclusion}

Bennett 3R chains have appeared in the literature~\cite{Selig_2014,Hegedus_2015} and some of their properties are known: Their configuration variety spans a 5-dimensional projective space and is the complete intersection of this space with two quadrics. In the context of bond theory and overconstrained closed kinematic 6R chains~\cite{Hegedus_2013b}, they give rise to ``short bonds''.

We complemented these results with several new insights: The characterization of kinematically equivalent Bennett 5-spaces in terms of their kinematic invariant \( K \) shows that different Bennett 3R chains can have the same Bennett 5-space up to kinematic equivalence. We also showed that Bennett 3R chains with the same ordered pair of distinct oriented outer axes have exactly the same Bennett 5-space, whatever their middle axis, and that such middle axes exist in abundance. A description of all possible middle axes for given outer axes is the subject of future work. Bennett 5-spaces are special in the sense that the back-projection of their polar line describes just a rotation (and not a ``true'' vertical Darboux motion) in the kinematic model of \( \bbD\bbH \) as a projective space over the dual numbers. Up to kinematic equivalence, this property characterizes them. Moreover, we provided simple equations for two quadrics that determine the configuration variety within a Bennett 5-space.

While we believe that this study is interesting in its own right, we intend to use it for the \( k \)-pose synthesis of overconstrained closed kinematic 6R chains. These are linkages of six revolute joints arranged in a loop. They can move, although generic linkages of this type are rigid. Synthesis in this context amounts to specifying the six axes in such a way that the configuration curve on the Study quadric passes through \( k \) prescribed points. This is considered to be a difficult problem of mechanism science for which few explicit constructions are known. Mechanisms obtained in this way have recently been shown to have some advantageous features~\cite{Huczala_2024,Huczala_2025,Postulka_2025}.

\section*{Funding}

This research was funded in whole or in part by the Austrian Science Fund (FWF) [10.55776/I6233]. For open access purposes, the authors have applied a CC BY public copyright license to any author-accepted manuscript version arising from this submission.

Zijia Li is supported by the National Key R\&D Program of China [2023YFA1009401] and partially supported by the Open Foundation of the State Key Laboratory of Mathematical Sciences [SLMS-2025-KFKT-TD-08], the Strategic Priority Research Program of the Chinese Academy of Sciences [XDA0480501] and the Guangdong Basic and Applied Basic Research Foundation [2024A1515010506].

The funders had no role in the study design, the analysis, the writing of the article, or the decision to submit it for publication.

\section*{\raggedright{}Declaration of generative AI and AI-assisted technologies in the manuscript preparation process}

During the preparation of this work the authors used Claude (Anthropic) to write the code for the figures and to suggest improvements to wording and language style. After using this tool, the authors reviewed and edited the content as needed and take full responsibility for the content of the published article.

\printbibliography{}

\end{document}